\documentclass[a4paper,11pt,reqno,twoside]{amsart}
\usepackage{amsmath}
\usepackage{amsfonts}
\usepackage{amssymb}
\usepackage{amsthm}
\usepackage{verbatim}
\usepackage{hyperref}
\theoremstyle{plain}
\newtheorem{thm}{Theorem}[section]
\newtheorem{cor}[thm]{Corollary}
\newtheorem{lemma}[thm]{Lemma}
\newtheorem{prop}[thm]{Proposition}

\theoremstyle{definition}
\newtheorem{remark}[thm]{Remark}

\numberwithin{equation}{section}

\DeclareMathOperator{\supp}{supp}
\DeclareMathOperator{\res}{Res}
\newcommand{\R}{\mathbb{R}}
\newcommand{\N}{\mathbb{N}}
\newcommand{\C}{\mathbb{C}}
\newcommand{\Z}{\mathbb{Z}}

\newcommand{\al}{\alpha}

\newcommand{\de}{\delta}
\newcommand{\si}{\sigma}
\newcommand{\la}{\lambda}

\newcommand{\Om}{\Omega}

\newcommand{\ep}{\varepsilon}
\newcommand{\vp}{\varphi}

\newcommand{\Pol}{\mathcal{P}}

\newcommand{\uint}{\int_0^\infty}
\newcommand{\ra}{\rightarrow}
\newcommand{\da}{\downarrow}
\newcommand{\bigO}{\mathcal{O}}
\newcommand{\qbin}[2]{\genfrac{[}{]}{0pt}{}{#1}{#2}}
\newcommand{\rps}[5]{\, {}_{#1}\vp_{#2}\left(
\genfrac{.}{.}{0pt}{}{#3}{#4};q, {#5}\right)}

\newcommand{\Ker}{\text{\rm Ker}}

\title[Self-adjoint difference operators]{Self-adjoint difference
  operators and classical solutions to the Stieltjes--Wigert moment
  problem}

\author{Jacob S.~Christiansen and Erik Koelink}
\address{Katholieke Universiteit Leuven, Departement Wiskunde,
Celestijnenlaan 200B, B-3001 Leuven, Belgium}
\email{stordal@wis.kuleuven.ac.be}
\address{Technische Universiteit Delft, DIAM, PO Box 5031,
2600 GA Delft, the Netherlands}
\email{h.t.koelink@ewi.tudelft.nl}

\begin{document}

\begin{abstract}
The Stieltjes--Wigert polynomials, which correspond to an
indeterminate moment problem on the positive half-line, are
eigenfunctions of a second order $q$-difference operator. We consider
the orthogonality measures for which the difference operator is
symmetric in the corresponding weighted $L^2$-spaces.
Under some additional assumptions these measures
are exactly the solutions to the $q$-Pearson equation.
In the case of discrete and absolutely continuous measures the
difference operator is essentially self-adjoint, and the corresponding
spectral decomposition is given explicitly. In particular, we find an
orthogonal set of $q$-Bessel functions complementing the
Stieltjes--Wigert polynomials to an orthogonal basis for $L^2(\mu)$
when $\mu$ is a discrete orthogonality measure solving the $q$-Pearson
equation. To obtain the spectral decomposition of the difference
operator in case of an absolutely continuous orthogonality measure we
use the results from the discrete case combined with direct integral
techniques.
\end{abstract}

\maketitle

\begin{small}
\end{small}

\begin{small}
\end{small}

\section{Introduction}

As part of the Askey-scheme \cite{K&S} of basic hypergeometric orthogonal
polynomials, the Stieltjes--Wigert polynomials are eigenfunctions of a
second-order $q$-difference operator. This operator is given by
\[
\bigl(Lf\bigr)(x)=f(xq)-\frac{1}{x}\,f(x)+\frac{1}{x}\,f(x/q)
\]
or, in a more compact form,
\[
L=T_q-x^{-1}(I-T_{q^{-1}}),
\]
where $T_a$ denotes the operator defined by $\bigl(T_a f\bigr)(x)
=f(ax)$ for fixed $a\neq 0$. We always take $q$ as a fixed number in
$(0,1)$. Clearly, $L$ preserves the space of polynomials.

In this paper we consider $L$ as a (possibly) unbounded operator on $L^2(\mu)$,
where $\mu$ is assumed to be a solution to the Stieltjes--Wigert
moment problem, i.e. a positive measure on $[0,\infty)$ such that
\begin{equation}\label{eq:SWmoments}
\int_0^\infty x^n d\mu(x)=q^{-\binom{n+1}{2}}, \quad n\geq 0.
\end{equation}
Since the Stieltjes--Wigert moment problem is indeterminate, there
are infinitely many positive measures to choose from. The operator
$(L,\mathcal{P})$ with domain the space $\Pol$ of polynomials is
always symmetric on $L^2(\mu)$. However, the polynomials are only
dense in $L^2(\mu)$ when $\mu$ is a so-called $N$-extremal solution
to the moment problem, see e.g. \cite[Chapter 2]{Akhi}. So instead
we consider $L$ with a larger domain $L(D)$ which will be specified
in \eqref{eq:defdomainL}. Under certain restrictions on $T_{q^{\pm
1}}$, this operator turns out only to be symmetric for a special
class of solutions to the moment problem, namely the solutions that
satisfy the $q$-Pearson equation or, in the setup of
\cite{ChriJMAA}, the solutions that are fixed points of the
transformation $T$ defined in \cite[Def.~2.4]{ChriJMAA}. Such
solutions are also called ``classical'' in \cite{ChriJMAA}. We give
the precise condition that $\mu$ has to satisfy in Proposition
\ref{prop:tau}.

The question now raises if $L$ can be extended to a self-adjoint
operator on $L^2(\mu)$ when $\mu$ is a classical solution to the
moment problem. We deal with the cases of discrete solutions,
respectively absolutely continuous solutions, in Section
\ref{discrete} and Section \ref{continuous}.

In Section \ref{discrete}, where $\mu$ is supposed to be discrete, we
show that $L$ is unitarily equivalent to a doubly infinite Jacobi
operator acting on $\ell^2(\Z)$. The theory of unbounded Jacobi
operators then leads to the fact that $L$ is essentially self-adjoint.
Starting from two explicit eigenfunctions of $L$ constructed in
Section 2, the spectrum of $L$ is computed in Theorem
\ref{thm:spectrum}. The spectrum is purely discrete (except for the
point $0$) and has an unbounded negative part and a bounded positive
part. The positive part is simple and each point corresponds to a
Stieltjes--Wigert polynomial of fixed degree. The negative part is
also simple and each point corresponds now to a $q$-Bessel function of
the second kind. This leads to orthogonality relations for the
Stieltjes--Wigert polynomials and for Jackson's second $q$-Bessel
functions. None of the discrete measures under consideration are
canonical solutions in the sense of \cite[Def.~3.4.2, p.~115]{Akhi},
and hence the space of polynomials has codimension $+\infty$ in the
corresponding weighted $L^2$-spaces. Our analysis leads to an explicit
set of orthogonal functions complementing the Stieltjes--Wigert
polynomials to a basis for $L^2(\mu)$.

In the case where $\mu$ is absolutely continuous, the operator $L$ is
again essentially self-adjoint. We show this in Section
\ref{continuous} using direct integrals of Hilbert spaces and the
results of Section \ref{discrete}. The spectrum of $L$ has a purely
discrete positive part, where each point is of infinite multiplicity
and corresponds to a Stieltjes--Wigert polynomials of fixed degree
times an arbitrary $q$-periodic function, i.e. a function $f$
satisfying $f(xq)=f(x)$ for all $x>0$. In case
$\text{supp}(\mu)=[0,\infty)$, the continuous spectrum of $L$ is
$(-\infty, 0]$ and each point here is simple. We also give an explicit
formula for the spectral measure. The approach in Section
\ref{continuous} should be compared with related ideas of Berg
\cite{Berg}.

The indeterminate cases within the Askey-scheme have been classified
in \cite{ChriPhD} and one may ask if a similar construction is
possible for other cases as well. For the $q$-Laguerre polynomials the
analysis is already done in \cite{CiccKK}, where the motivation comes
from quantum groups and limit transitions of the big $q$-Jacobi
polynomials. Formal limit results of \cite{CiccKK} lead to the results
of Section \ref{discrete}, and we note that the methods of Section
\ref{continuous} can be used for the $q$-Laguerre case as well. See
also \cite{ChriCA} for the transformation corresponding to the
$q$-Pearson equation. For other cases in the indeterminate part of the
Askey-scheme several problems arise, and it is not clear if symmetry
of the difference operator for the corresponding orthogonal
polynomials has a clear-cut meaning for solutions to the moment
problem.

{\bf Acknowledgement.} We thank the referee for useful
suggestions, and Barry Simon for a remark that led to an improvement
of Section \ref{sec:diffoperator}.

\section{Difference operator}\label{sec:diffoperator}

\subsection{Difference operator}
Consider the second order $q$-difference operator
\begin{equation}
\label{eq:diffopSW} \bigl( Lf\bigr) (x) = f(xq) - \frac{1}{x} \,
f(x) + \frac{1}{x} \, f(x/q).
\end{equation}
The motivation for studying $L$ is the fact that the Stieltjes--Wigert
polynomials
\begin{equation}
\label{SW}
S_n(x;q)=\frac{1}{(q;q)_n}\sum_{k=0}^n\genfrac{[}{]}{0pt}{}{n}{k}_q
(-1)^kq^{k^2}x^k, \quad n=0,1,\ldots
\end{equation}
are eigenfunctions of $L$ corresponding to the eigenvalues $q^n$,
see Proposition \ref{prop:eigenfunctions} below. Here we use the
notation
\[
(q;q)_0=1, \quad (q;q)_n=\prod_{k=1}^n(1-q^k), \quad n=1,2,\ldots
\]
and
\[
\genfrac{[}{]}{0pt}{}{n}{k}_q=\frac{(q;q)_n}{(q;q)_k(q;q)_{n-k}},
\quad 0\leq k\leq n.
\]
Throughout the paper we assume that $0<q<1$ and follow the notation
of Gasper and Rahman \cite{GaspR} for basic hypergeometric series.

Recall that the image measure $\tau(\mu)$ of a finite positive measure
$\mu$ under a measurable map $\tau$ is defined by
\[
\tau(\mu)(A)=\mu\bigl(\tau^{-1}(A)\bigr)
\]
for any measurable set $A$. Recall also that integration with respect
to $\tau(\mu)$ is carried out via the rule
\[
\int f\,d\tau(\mu)=\int (f\circ\tau)\,d\mu.
\]
In what follows we denote by $\tau_a:(0,\infty)\ra (0,\infty)$ the map
given by $x\mapsto ax$ for fixed $a>0$.

Writing $M$ for the operator of multiplication by ${1}/{x}$, we
see that $L$ can be written as 
\[
L=T_q - M + M\circ T_{q^{-1}}.
\]
Our first task is therefore to define and discuss the operators $M$
and $T_{q^{\pm 1}}$ as possibly unbounded operators on $L^2(\mu)$,
where $\mu$ for the time being is supposed to be any finite positive
(Borel) measure on $(0,\infty)$. We define the operator $M$ on the
maximal domain
\[
D(M)=\Bigl\{f\in L^2(\mu)\,\Big|\, 
\uint \frac{1}{x^2} | f(x)|^2\, d\mu(x) <\infty\Bigr\}.
\]
As regards the operators $T_{q^{\pm 1}}$, it may happen that one (or
both) of them is identically zero on $L^2(\mu)$. This happens if $xq$
(or $x/q$) never belongs to $\supp(\mu)$ when $x\in\supp(\mu)$ (and
hence for example if $\mu$ is discrete and supported on $\{tq^{2n}\mid
n\in\Z \}$ for some $t>0$). To avoid this situation we require that
$T_{q^{\pm 1}}$, defined on the maximal domains
\[
D(T_{q^{\pm 1}})=\bigl\{f\in L^2(\mu)\mid T_{q^{\pm 1}}f\in L^2(\mu)\bigr\},
\]
have trivial kernels, i.e. $\Ker(T_{q^{\pm 1}})=\{0\}$. For any
Borel set $A\subset(0,\infty)$, the indicator function $\chi_A$
belongs to $D(T_{q^{\pm 1}})$ since
\[
\uint |(T_{q^{\pm 1}}\chi_A)(x)|^2\,d\mu(x)=\mu(q^{\mp 1}A)<\infty.
\]
When $\mu(A)>0$, we have $\chi_A\not=0$ in $L^2(\mu)$ and the
requirement on the kernels therefore implies that $\mu(q^{\mp
  1}A)=\tau_{q^{\pm 1}}(\mu)(A)>0$. In other words, $\mu$ is
absolutely continuous with respect to $\tau_{q^{\pm 1}}(\mu)$, that
is, $\tau_{q^{\pm 1}}$ preserve the support of $\mu$. Note that the
domains $D(T_{q^{\pm 1}})$ are dense in $L^2(\mu)$ since the set of
finite linear combinations of indicator functions is dense in
$L^2(\mu)$.

With the above assumptions in mind we define $L$ as the possibly
unbounded operator on $L^2(\mu)$ with domain
\begin{equation}\label{eq:defdomainL}
D(L)=\bigl\{ f\in L^2(\mu) \mid f\in D(T_q)\cap D(M)\cap D(T_{q^{-1}}), \,
T_{q^{-1}}f \in D(M) \bigr\}.
\end{equation}

\begin{prop}
\label{prop:tau} Let $\mu$ be a positive measure on $(0,\infty)$
such that 
\[
m_n:=\uint x^n\, d\mu(x) <\infty \quad \mbox{for } n\geq -2.
\]
Assume that $T_{q^{\pm 1}}\colon D(T_{q^{\pm 1}})\to L^2(\mu)$ have
trivial kernels. Then the domain $D(L)$ defined in
\eqref{eq:defdomainL} is dense in $L^2(\mu)$ and the operator $(L,
D(L))$ is symmetric on $L^2(\mu)$ if and only if the measure
$\tau_q(\mu)$ is absolutely continuous with respect to $\mu$ and the
Radon--Nikodym derivative is given by
\begin{equation}
\label{eq:tau}
\frac{d\tau_q(\mu)}{d\mu}=\frac{1}{x} \quad\mbox{a.e. with respect to } \mu.
\end{equation}
\end{prop}

\begin{remark}
\label{rem:moments}
When $\mu$ is a finite positive measure on $(0,\infty)$ satisfying
\eqref{eq:tau}, it follows by induction that $\tau_{q^n}(\mu)$ is
absolutely continuous with respect to $\mu$ for all $n\in\Z$ and
\[
\frac{d\tau_{q^n}(\mu)}{d\mu}=\frac{q^{\binom{n}{2}}}{x^n}
\quad\mbox{a.e. with respect to } \mu.
\]
This in particular means that $\mu$ has moments of all orders and if
$\mu$ is a probability measure, then
\[
\uint x^n\,d\mu(x)=q^{-\binom{n+1}{2}} \quad\mbox{for all } n\in\Z.
\]
So the requirement in Proposition \ref{prop:tau} on the existence
of the first two negative moments is actually implied by \eqref{eq:tau}.
Moreover, we see that $\mu$ is uniquely determined by its restriction
$\mu|_{(q,1]}$ to the interval $(q,1]$ (or any other
interval of the form $(tq^{k+1},tq^k]$ for $t>0$ and $k\in\Z$). See
\cite[Section 2]{ChriJMAA} for more details.
\end{remark}

\begin{proof}
Since by assumption $m_{-2}<\infty$, we see that $\chi_A\in D(M)$ for
any Borel set $A\subset(0,\infty)$. We have already observed that
$\chi_A\in D(T_{q^{\pm 1}})$ and that $T_{q^{-1}}\chi_A = \chi_{qA}\in
D(M)$.  Hence, all indicator functions are contained in $D(L)$, and
finite linear combinations of these functions are dense in $L^2(\mu)$.

Suppose that $f, g\in D(L)$, then
\begin{align*}
\langle Lf,g\rangle &=
\uint\bigl(Lf\bigr)(x)\,\overline{g(x)}\,d\mu(x) \\
&=\uint\Bigl(f(xq)-\frac{1}{x}\,f(x)+\frac{1}{x}\,f(x/q)\Bigr)\,
\overline{g(x)}\,d\mu(x) \\
&=\uint f(x)\,\overline{g(x/q)}\,d\tau_q(\mu)(x)
-\uint f(x)\,\frac{\overline{g(x)}}{x}\,d\mu(x)
+\uint f(x)\,\frac{\overline{g(xq)}}{xq}\,d\tau_{q^{-1}}(\mu)(x),
\end{align*}
using the fact that each term is integrable.
The right-hand side can be written as $\langle f,Lg\rangle$ if
and only if
\begin{equation}\label{eq:symmetry}
\begin{split}
\uint f(x) \overline{g(qx)}\, d\mu(x)
+& \uint f(x) \overline{\frac{g(x/q)}{x}}\, d\mu(x) = \\
 &\uint f(x)\,\overline{g(x/q)}\,d\tau_q(\mu)(x)
+\uint f(x)\,\frac{\overline{g(xq)}}{xq}\,d\tau_{q^{-1}}(\mu)(x).
\end{split}
\end{equation}
Now, if $\tau_q(\mu)$ and $\tau_{q^{-1}}(\mu)$ are both absolutely
continuous with respect to $\mu$ and the conditions
\[
\frac{d\tau_q(\mu)}{d\mu}=\frac{1}{x} \quad\mbox{and}\quad
\frac{d\tau_{q^{-1}}(\mu)}{d\mu}=xq
\quad\mbox{a.e. with respect to } \mu
\]
are met, then \eqref{eq:symmetry} is satisfied. Since
$\tau_{q^{-1}}=\tau_q^{-1}$, these conditions are equivalent and the
``if'' part of the proposition follows.

Conversely, if $(L,D(L))$ is symmetric, then \eqref{eq:symmetry} holds
for all $f,g\in D(L)$. Take $f=\chi_A$, $g=\chi_B$, then
\[
\int_{A\cap q^{-1}B}  \, d\mu(x)
+ \int_{A\cap qB} \frac{1}{x}\, d\mu(x) =
\int_{A\cap qB} \,d\tau_q(\mu)(x)
+\int_{A\cap q^{-1}B} \frac{1}{xq}\,d\tau_{q^{-1}}(\mu)(x).
\]
Now take $A\subset (q^{k+1},q^k]$ for some $k\in\Z$, and
set $B=q^{-1}A$ or $A=qB$. This gives $A\cap q^{-1}B=\emptyset$
and therefore
\[
\int_{A} \frac{1}{x}\, d\mu(x) =
\tau_q(\mu)(A).
\]
Since any Borel set $A\subset(0,\infty)$ can be written as a
disjoint union $A=\cup_{k\in\Z} A_k$, where $A_k=A \cap
(q^{k+1},q^k]$, we find that
\[
\tau_q(\mu)(A) = \sum_{k\in\Z} \tau_q(\mu)(A_k)
=\sum_{k\in\Z} \int_{A_k} \frac{1}{x}\, d\mu(x)
= \int_A \frac{1}{x}\, d\mu(x),
\]
recalling that $1/x$ is integrable with respect to $\mu$. In particular,
$\tau_q(\mu)$ is absolutely continuous with respect to $\mu$ and
\eqref{eq:tau} is satisfied.
\end{proof}

\begin{remark}
  When $\mu$ is an $N$-extremal (or $m$-canonical) solution to the
  Stieltjes--Wigert moment problem, then $\tau_{q^{\pm 1}}$ do not
  preserve the support of $\mu$. See \cite[Section 3]{ChriJMAA} for
  details. So the assumptions on $T_{q^{\pm 1}}$ in Proposition
  \ref{prop:tau} exclude canonical solutions of all orders.
\end{remark}

In this paper we shall mainly focus on discrete and absolutely
continuous measures and state therefore the following consequence of
Proposition \ref{prop:tau}. As for notation, we denote by $\delta_x$
the unit mass at the point $x$.

\begin{cor}
\label{cor:symmetricmeasure}
\par\noindent
{\rm (i)} Suppose that $t>0$ and let $\mu_t$ be a positive discrete
measure of the form
\[
\mu_t=\sum_{k=-\infty}^\infty m_t(k) \de_{tq^k},
\]
where $m_t(k)>0$ for all $k\in\Z$ and $\sum_{k=-\infty}^\infty
m_t(k)<\infty$. The operator $L$ is symmetric on $L^2(\mu_t)$ if and
only if
\begin{equation}
\label{eq:discrete}
m_t(k+1)=tq^{k+1}m_t(k) \mbox{ for all } k\in\Z.
\end{equation}
\par\noindent
{\rm (ii)} Let $\mu$ be an absolutely continuous measure on
$(0,\infty)$ given by a positive density function $w$ satisfying
$\uint w(x)dx<\infty$. Assume that $\mu$ and $\tau_{q^{\pm 1}}(\mu)$
have the same support. The operator $L$ is symmetric on $L^2(\mu)$
if and only if
\begin{equation}
\label{eq:qPearson}
w(xq)=xw(x) \text{ for all } x\in(0,\infty).
\end{equation}
\end{cor}

\begin{remark}
\label{rmk:SWmoments}
(i) The condition \eqref{eq:discrete} is equivalent to
$m_t(k)=t^kq^{\binom{k+1}{2}}m_t(0)$ for $k\in\Z$.
If we set $1/m_t(0)=(-tq,-1/t,q;q)_\infty$, it follows by the
triple product identity \cite[(1.6.1)]{GaspR} that $\mu_t$ becomes a
probability measure.
\par\noindent
(ii) The condition \eqref{eq:qPearson} is the $q$-Pearson equation for
the Stieltjes--Wigert polynomials, see e.g. \cite{Marcellan+Medem} and
\cite{Nodarse+Medem}. This equation is for example satisfied by the
log-normal density
\[
w(x)=\frac{1}{\sqrt{x}}e^{\tfrac{1}{2}\tfrac{(\log x)^2}{\log q}},
\quad x>0
\]
and (for fixed $c>0$) by the infinite products
\[
w_c(x)=\frac{x^{c-1}}{(-q^{1-c}x,-q^c/x;q)_\infty}, \quad x>0.
\]
Note also that \eqref{eq:qPearson} is invariant under multiplication
with $q$-periodic functions, that is, functions which satisfy
$f(xq)=f(x)$ for $x>0$.
\end{remark}

In the setting of Proposition \ref{prop:tau} we find
\[
\int_0^\infty |f(xq)|^2\,d\mu(x)=
\int_0^\infty\frac{1}{x}|f(x)|^2\,d\mu(x)=
q\int_0^\infty\frac{1}{x^2}|f(x/q)|^2\,d\mu(x),
\]
showing that $L$ is well-defined on any continuous function $f$
satisfying $f(x)=\bigO(x^N)$ as $x\ra\infty$ and $f(x)=\bigO(x^{-M})$
as $x\ra 0$ for some $N, M\geq 0$, cf. Remark \ref{rem:moments}.

\subsection{Eigenfunctions}

The ${}_1\vp_1$-series with lower parameter equal to zero, say
$\rps{1}{1}{a}{0}{y}$, satisfies the second order $q$-difference
equation
\begin{equation}\label{eq:difffor1vp1}
-ay\, f(yq) + (y-q)\, f(y) + q\, f(y/q) = 0.
\end{equation}
This result can be obtained from the second order $q$-difference
equation for the ${}_2\vp_1$-series \cite[Exerc.~1.13]{GaspR} by
taking a limit.

By looking for solutions of the form
$\displaystyle{\sum_{k=0}^\infty c_k y^{\la+k}}$, respectively
$\displaystyle{\sum_{k=0}^\infty c_k y^{\la-k}}$, with $c_0=1$, we
see that
\begin{equation}\label{eq:1phi1solutions}
\rps{1}{1}{a}{0}{y} \quad\mbox{and}\quad
y^\al\rps{1}{1}{a}{0}{\frac{q^2}{y}} , \quad q^\al a=1
\end{equation}
both satisfy \eqref{eq:difffor1vp1}.

\begin{prop}\label{prop:eigenfunctions}
The functions defined by
\begin{equation*}
\phi_z(x) = \rps{1}{1}{1/z}{0}{-xzq},
\qquad \Phi_z(x) = x^{\ln z/\ln q}
\rps{1}{1}{1/z}{0}{-\frac{q}{xz}}
\end{equation*}
are solutions to the eigenvalue equation $Lf = zf$. Here
$\phi_z(x)$ is defined for $x,z\in \C$, where
the case $z=0$ has to be interpreted as the limit
\[
\phi_0(x)=\rps{0}{1}{-}{0}{-xq},
\]
and $\Phi_z(x)$ is defined for $x\in (0,\infty)$ and $z\in
\C\backslash (-\infty, 0]$.

In particular, the Stieltjes--Wigert polynomials are solutions to
the eigenvalue equations
\[
L S_n(\,\cdot\,;q)=q^n S_n(\,\cdot\,;q), \quad n=0,1,\ldots.
\]
\end{prop}

\begin{remark}
The function
\[
\phi_0(x)=\sum_{n=0}^\infty\frac{(-1)^nq^{n^2}x^n}{(q;q)_n}, \quad x\in\C
\]
is also known as the entire Rogers--Ramanujan function, since its
values at $-1$ and $-q$ appear in the celebrated identities
\cite[(2.7.3/4)]{GaspR}
\[
\sum_{n=0}^\infty\frac{q^{n^2}}{(q;q)_n}=\frac{1}{(q,q^4;q^5)_\infty}
\quad \mbox{and} \quad
\sum_{n=0}^\infty\frac{q^{n(n+1)}}{(q;q)_n}=\frac{1}{(q^2,q^3;q^5)_\infty}.
\]
The reader is referred to \cite{Andrews} and \cite{Hayman} for
interesting results about the zeros of $\phi_0$, which are all
positive and simple.
\end{remark}

\begin{proof}
The result follows from \eqref{eq:difffor1vp1} and
\eqref{eq:1phi1solutions} if we replace $a$ by $1/z$ and $y$ by
$-xzq$.
Since
\[
\sum_{k=0}^n\genfrac{[}{]}{0pt}{}{n}{k}_q (-1)^kq^{k^2}x^k
=\rps{1}{1}{q^{-n}}{0}{-q^{n+1}x},
\]
the last assertion follows immediately from \eqref{SW}.
\end{proof}

To get hold of the behavior of $\Phi_z(x)$ as $x\da 0$, we need the
following result.
\begin{lemma}\label{lem:x=0}
As $x\da 0$, we have
\[
\rps{0}{1}{-}{-zq/x}{-\frac{z^2q}{x}} \,\longrightarrow \,
\rps{0}{0}{-}{-}{z}=(z;q)_\infty,
\]
and the convergence is uniform for $z$ in compact subsets of
$\C\setminus(-\infty,0)$.
\end{lemma}
\begin{proof}
Notice that
\[
\rps{0}{1}{-}{-zq/x}{-\frac{z^2q}{x}}=1+\sum_{n=1}^\infty
\frac{(-1)^nq^{n^2}}{(q;q)_n}\frac{z^{2n}}{(x+zq)\cdots(x+zq^n)}
\]
for $z\in\C\setminus(-\infty,0)$ and $x>0$. The termwise
convergence is thus obvious. Let $K$ be a compact subset of
$\C\setminus(-\infty,0)$ and take $\de>0$ such that
$|z-t|\geq\de$ for all $z\in K$ and $t<0$. Clearly,
\[
|(x+zq)\cdots(x+zq^n)|\geq\de^nq^{\binom{n+1}{2}}
\]
and since the right-hand side is independent of $z\in K$ and
$x>0$, we have dominated convergence.
\end{proof}

A limit case of Heine's transformation formula for the
${}_2\vp_1$-series \cite[(0.6.8/9)]{K&S} tells us that
\begin{equation}
\label{eq:Heine}
\rps{1}{1}{1/z}{0}{-\frac{q}{xz}}=(-q/xz;q)_\infty
\rps{0}{1}{-}{-q/xz}{-\frac{q}{xz^2}}
\end{equation}
and according to Lemma \ref{lem:x=0}, the ${}_0\vp_1$-series on the
right-hand side converges to $(1/z;q)_\infty$ as $x\da 0$. We follow
the convention that in a fraction the part to the right of $/$ is the
denominator. So in (2.8), for example, we write $(-q/xz;q)_\infty$
instead of $(-\frac{q}{xz};q)_\infty$. The infinite product
$(-q/xz;q)_\infty$ does not have a limit as $x\ra 0$, but for $x=tq^n$
we have
\begin{equation}
\label{eq:inf-prod}
(-q/xz;q)_\infty=(-q^{1-n}/tz;q)_\infty
=\frac{(-tz;q)_n(-q/tz;q)_\infty}{(tz)^n q^{\binom{n}{2}}}.
\end{equation}

\section{Spectral analysis for the discrete case}
\label{discrete}
In this section we consider $L$ as an unbounded symmetric operator on
the Hilbert space $L^2(\mu_t)$, where $\mu_t$ is the discrete measure
from Corollary \ref{cor:symmetricmeasure} (i). Throughout the section
the parameter $t>0$ will be fixed.
\subsection{$\ell^2(\Z)$ setup}
\label{ell2}
Since $L^2(\mu_t)$ essentially is a weighted $\ell^2$-space over the
integers, we start by defining a unitary operator $U\colon
L^2(\mu_t)\to \ell^2(\Z)$ by
\[
Uf = \sum_{k=-\infty}^\infty f(tq^k) \sqrt{m_t(k)}\, e_k,
\]
where $\{ e_k\}_{k\in\Z}$ denotes the standard orthonormal
basis for $\ell^2(\Z)$. The adjoint of $U$ is given by
\[
\bigl(U^\ast e_k\bigr)(tq^r)=\frac{1}{\sqrt{m_t(k)}}\de_{k,r}
\]
and the operator $J=ULU^\ast$ becomes a doubly infinite Jacobi
operator on $\ell^2(\Z)$. More precisely, $J$ has the form
\[
Je_k=a_k e_{k+1}+b_k e_k+a_{k-1} e_{k-1}, \quad k\in\Z
\]
with
\[
a_k=
\frac{1}{\sqrt{tq^{k+1}}}
\quad\mbox{and}\quad b_k = -\frac{1}{tq^k}.
\]
In what follows, we denote by $\mathcal{D}$ the subspace of
$\ell^2(\Z)$ consisting of finite linear combinations of the basis
elements. Clearly, $(J, \mathcal{D})$ is a densely defined symmetric
operator on $\ell^2(\Z)$. But more importantly, we have the following
result.
\begin{thm}
\label{thm:self-adjoint}
The operator $(J, {\mathcal D})$ is essentially self-adjoint.
\end{thm}

By the unitary intertwiner $U$, the operator $(J, {\mathcal D})$
corresponds to $(L, U^\ast{\mathcal D}U)$ which is a restriction of
the operator $(L,D(L))$ considered in Proposition \ref{prop:tau}. The
domain $U^\ast{\mathcal D}U$ consists of the compactly supported
functions in $L^2(\mu)$, and it is straightforward to check that this
is a core for the closure of $(L,D(L))$. So by the above theorem,
$(L,D(L))$ is essentially self-adjoint in the case $\mu=\mu_t$.

\begin{proof} We employ a theorem of Masson and Repka
\cite{MassR}, see also \cite[Thm. 4.2.2]{KoelLaredo}. For this we
define the operators
\[
J^\pm := P^\pm J\big\vert_{{\mathcal D}^\pm},
\]
where $P^+$ and $P^-$ are the orthogonal projections onto
$\text{span}\{ e_k\mid k\geq 0\}$, respectively $\text{span}\{
e_k\mid k< 0\}$, and
\[
{\mathcal D}^+= {\mathcal D}\cap \text{span}\{ e_k\mid k\geq 0\},
\quad{\mathcal D}^-= {\mathcal D}\cap \text{span}\{ e_k\mid k<0\}.
\]
Notice that $J^\pm$ are Jacobi operators on $\ell^2(\N)$ with finite
linear combinations of the basis vectors as domain. The theorem of
Masson and Repka states that the deficiency indices of $J$ can be
obtained by adding the deficiency indices of $J^+$ and $J^-$, see e.g.
Akhiezer \cite[Ch.~4]{Akhi} or Berezanski\u\i\ \cite[Ch.~7]{Bere} for
more information. The deficiency indices of $J^-$ are $(0,0)$ since
the coefficients $a_k$ and $b_k$ are bounded as $k\to-\infty$. For the
deficiency indices of $J^+$ we observe that $a_k+b_k+a_{k-1}$ is
bounded from above
for $k\geq 0$, and by \cite[Addenda and problems to Chap.~1]{Akhi} or
\cite[Thm. 1.4, p.~505]{Bere} this implies that $J^+$ is essentially
self-adjoint.  Hence, the deficiency indices of $J^+$ are $(0,0)$ and
we conclude that the deficiency indices of $J$ are also $(0,0)$. The
statement follows.
\end{proof}

The closure of $(J,{\mathcal D})$ thus coincides with the adjoint
operator $(J^\ast,{\mathcal D}^\ast)$, which is defined on the maximal
domain
\[
{\mathcal D}^\ast = \Bigl\{ v\in\ell^2(\Z): \sum_{k=-\infty}^\infty
  \bigl|a_k v_{k+1} + b_k v_k + a_{k-1}v_{k-1}\bigr|^2 <\infty \Bigr\}.
\]
\subsection{Wronskian and Green function}
\label{WandG}
We now aim at finding the spectrum of the self-adjoint operator
$(J^\ast, \mathcal{D}^\ast)$. In this connection the functions from
Proposition \ref{prop:eigenfunctions} become very useful. We set
\[
\psi_k(z)=t^{k/2}q^{k(k+1)/4}\phi_z(tq^k),
\]
respectively
\[
\Psi_k(z)=t^{k/2}q^{k(k+1)/4}\Phi_z(tq^k)/t^{\ln z/\ln q},
\]
and consider the two sequences $\psi(z)=\{\psi_k(z)\}_{k\in\Z}$ and
$\Psi(z)=\{\Psi_k(z)\}_{k\in\Z}$. Notice that $\psi(z)$ belongs to
$\ell^2$ as $k\ra\infty$ for all $z\in\C$, whereas $\Psi(z)$ belongs
to $\ell^2$ as $k\ra -\infty$ for $z\in\C\setminus\{0\}$. However,
except for special values to be determined later on, neither $\psi(z)$
nor $\Psi(z)$ is an element of $\ell^2(\Z)$. Since we divide by
$t^{\ln z/\ln q}$ in the definition of $\Psi_k(z)$, the sequence
$\Psi(z)$ is well-defined for all $z\in\C\setminus \{0\}$.

It follows from Proposition \ref{prop:eigenfunctions} that
$\psi(z)$ and $\Psi(z)$ are solutions to the eigenvalue equation
$Jv=zv$. Their Wronskian, i.e. the sequence defined by
\begin{equation}
\label{eq:Wronskian}
[\psi(z),\Psi(z)]_k=
a_k\bigl(\psi_{k+1}(z)\Psi_k(z)-\psi_k(z)\Psi_{k+1}(z)\bigr),
\quad k\in\Z,
\end{equation}
is therefore independent of $k$.

\begin{lemma}
\label{Lem:Wronskian}
The Wronskian of $\psi(z)$ and $\Psi(z)$ is given by
\[
[\psi(z),\Psi(z)]=-z(-tzq,-1/tz,1/z;q)_\infty.
\]
\end{lemma}
\begin{proof}
Inserting the expressions for $a_k$, $\psi_k(z)$ and $\Psi_k(z)$
in \eqref{eq:Wronskian}, we get after a few computations
\begin{align*}
[\psi(z),\Psi(z)]_k=z^{k}t^kq^{\binom{k+1}{2}}
\Biggl\{&\rps{1}{1}{1/z}{0}{-tzq^{k+2}}
\rps{1}{1}{1/z}{0}{-\frac{q^{1-k}}{tz}}\\
&-z\,\rps{1}{1}{1/z}{0}{-tzq^{k+1}}
\rps{1}{1}{1/z}{0}{-\frac{q^{-k}}{tz}}\Biggr\}.
\end{align*}
Since the Wronskian is independent of $k$, we evaluate the
expression by taking the limit $k\ra\infty$. Clearly, the
${}_1\varphi_1$-series with argument $-tzq^{k+2}$ (or $-tzq^{k+1}$)
converges to $1$ as $k\ra\infty$. Combining \eqref{eq:Heine} with
Lemma \ref{lem:x=0} and \eqref{eq:inf-prod}, we find that
\[
\rps{1}{1}{1/z}{0}{-\frac{q^{1-k}}{tz}} \sim
\frac{(-tz,-q/tz,1/z;q)_\infty}{(tz)^k q^{\binom{k}{2}}} \quad
\mbox{as } k\ra\infty,
\]
respectively
\[
\rps{1}{1}{1/z}{0}{-\frac{q^{-k}}{tz}} \sim
\frac{(-tz,-q/tz,1/z;q)_\infty}{(tz)^{k+1} q^{\binom{k+1}{2}}}
\quad \mbox{as } k\ra\infty,
\]
where ${\sim}$ means that the ratio of the right-hand side and the
left-hand side converges to $1$ as $k\ra\infty$. Therefore,
\[
[\psi(z),\Psi(z)]=\lim_{k\ra\infty}
\bigl(q^k-1/t\bigr) (-tz,-q/tz,1/z;q)_\infty
=-z(-tzq,-1/tz,1/z;q)_\infty
\]
and the desired result is established.
\end{proof}
With the Wronskian of $\psi(z)$ and $\Psi(z)$ at hand, we define the
Green function by
\[
G_z(j,l)=\frac{1}{[\psi(z),\Psi(z)]}
\begin{cases}
\psi_j(z)\Psi_l(z), \quad l\leq j, \\
\psi_l(z)\Psi_j(z), \quad l>j.
\end{cases}
\]
The resolvent of $(J^\ast,\mathcal{D}^\ast)$ is closely related to the
Green function, see e.g. \cite[Section 4.3]{KoelLaredo}. For any sequence
$v\in\ell^2(\Z)$, we have
\begin{equation}
\label{eq:resolvent} \bigl((J^\ast-z)^{-1}v\bigr)_j=
\sum_{l=-\infty}^\infty G_z(j,l)\,v_l, \quad z\in\C\setminus\R.
\end{equation}
\subsection{Spectral decomposition}
We denote by $E$ the resolution of the identity corresponding to the
self-adjoint operator $(J^\ast,\mathcal{D}^\ast)$. From general theory
(see e.g. \cite[Thm.~XII.2.10]{DunfS}) we know that
\begin{equation}
\label{eq:resolution} \bigl< E\bigl((a,b)\bigr)v,w\bigr>
=\lim_{\delta\da 0}\lim_{\ep\da 0} \frac{1}{2\pi
i}\int_{a+\delta}^{b-\delta} \bigl< (J^\ast-s-i\ep)^{-1}v,w \bigr>-
\bigl< (J^\ast-s+i\ep)^{-1}v,w \bigr>\, ds
\end{equation}
for $v, w\in\ell^2(\Z)$ and because of \eqref{eq:resolvent}, the inner
products in the integral can be written as
\begin{equation}
\label{eq:lj} \bigl< \bigl(J^\ast-(s\pm i\ep)\bigr)^{-1}v,w \bigr>=
\sum_{l\leq j}\frac{\psi_j(s\pm i\ep)\Psi_l(s\pm i\ep)} {[\psi(s\pm
i\ep),\Psi(s\pm i\ep)]}
(v_l\overline{w}_j+v_j\overline{w}_l)(1-\tfrac{1}{2}\delta_{j,l}).
\end{equation}
Since $\psi_k(z)$ is entire and $\Psi_k(z)$ is analytic in
$\C\setminus\{0\}$, it therefore follows that the spectral measure is
discrete and supported on the zeros of the Wronskian
$[\psi(z),\Psi(z)]$. We can read off these zeros from Lemma
\ref{Lem:Wronskian} and get $0$, $-q^r/t$ for $r\in\Z$ and $q^n$ for
$n\in\Z_+$.

\begin{thm}
\label{thm:spectrum}
The spectrum of $J^\ast$ is given by $\sigma(J^\ast)=-q^\Z/t\,\cup\,
\{0\}\,\cup\,q^{\Z_+}$. The accumulation point $0$ does not
belong to the point spectrum $\sigma_p(J^\ast)$.
\end{thm}
\begin{proof}
It is only left to prove that $0$ does not belong to the point
spectrum of $J^\ast$. We show that no non-trivial solution to the
equation $J v=0$ belongs to $\ell^2(\Z)$. In the end of the proof we
use the implication $\phi_0(t)=0 \Rightarrow\phi_0(tq)\neq 0$, which
follows from the fact that the zeros of $\phi_0$ are very well
separated, see e.g.  \cite[Section 3]{ChriJMAA}.

The space of solutions to the equation $a_k
v_{k+1}+b_k v_k+a_{k-1}v_{k-1}=0$ or, more explicitly,
\begin{equation}
\label{eq:nul}
v_{k+1}=\frac{1}{\sqrt{tq^{k-1}}}\,v_k-\sqrt{q}\,v_{k-1}, \quad k\in\Z
\end{equation}
is two-dimensional. We already know one solution, namely $\psi(0)$,
which is given by
\[
\psi_k(0)=t^{k/2}q^{k(k+1)/4}\phi_0(tq^k), \quad k\in\Z.
\]
Clearly $\psi(0)$ belongs to $\ell^2$ as $k\ra\infty$ but recalling
that $\phi_0(tq^{-2n})\sim (-1)^nt^nq^{-n^2}K(t)$ as $n\ra\infty$ for
some constant $K(t)>0$, see e.g. \cite{Mourad}, it follows that
\[
\psi_{-2n}(0)\sim (-1)^nq^{-n/2} K(t) \quad \mbox{as } n\ra\infty.
\]
Therefore, $\psi(0)$ does not belong to $\ell^2(\Z)$.

The sequence $\Psi(z)$ is not defined for $z=0$ so we need to look for
other solutions to \eqref{eq:nul}. Note that if $v_k$ has the form
\[
v_{k+1}=\frac{F_{k+1}}{t^{k/2}q^{k(k-1)/4}},
\]
then \eqref{eq:nul} is equivalent to
\[
F_{k+1}=F_k-tq^{k-1}F_{k-1}, \quad k\in\Z.
\]
With $F_0=0$ and $F_1=1$ (or, equivalently, $v_0=0$ and $v_1=1$) we
see that $F_k$, $k=0,1,\ldots$, essentially are $q$-Fibonacci polynomials
in $t$, see e.g. \cite{Carlitz}. In particular,
\[
F_{k+1}=\sum_{n=0}^{k-1} \qbin{k-n}{n}_q (-1)^nq^{n^2}t^n
\quad\mbox{and}\quad
F_k\ra\phi_0(t) \mbox{ as } k\ra\infty.
\]
There are two cases to be considered. 1) When $\phi_0(t)\neq 0$, the
solution to \eqref{eq:nul} with $v_0=0$ and $v_1=1$ does not belong to
$\ell^2$ as $k\ra\infty$.  Moreover, since this solution is linearly
independent of $\psi(0)$, there are no solutions to \eqref{eq:nul} in
$\ell^2(\Z)$. 2) In the case $\phi_0(t)=0$, the solution to \eqref{eq:nul}
with $v_0=0$ and $v_1=1$ is proportional to $\psi(0)$. But since
$\phi_0(tq)\neq 0$, the solution to \eqref{eq:nul} with $v_1=0$ and
$v_2=1$ is linearly independent of $\psi(0)$. This solution behaves
like $\phi_0(tq)/t^{k/2}q^{k(k-1)/4}$ as $k\ra\infty$ and as before we
see that no solution to \eqref{eq:nul} belongs to $\ell^2(\Z)$.
\end{proof}

\subsection{Orthogonality relations}
In this section we determine the spectral measure $E(\{\xi\})$ for
$\xi$ in the point spectrum of $J^\ast$. Our considerations will lead
to explicit orthogonality relations for the Stieltjes--Wigert
polynomials and the second $q$-Bessel functions of Jackson.

Along the way we will need the following auxiliary result.
\begin{lemma}
\label{lem:proportional}
For $c\in\C$ and $k, m\in\Z$, we have
\begin{equation}
\label{eq:prop}
(-c)^{m+k}\rps{1}{1}{-cq^{-m}}{0}{q^{1+m+k}}=
q^{m(m+k)}\rps{1}{1}{-cq^{-m}}{0}{q^{1-m-k}}.
\end{equation}
\end{lemma}
\begin{proof}
Because of symmetry it suffices to establish the identity for
$m+k\geq 0$. Applying the transformation \cite[(0.6.8/9)]{K&S}, we
see that the right-hand side of \eqref{eq:prop} can be written as
\[
q^{m(m+k)}\sum_{n=m+k}^\infty
\frac{(q^{1-m-k+n};q)_\infty}{(q;q)_n}(-c)^n q^{n(n-2m-k)}=
(-c)^{m+k}\sum_{n=0}^\infty\frac{(q^{1+m+k+n};q)_\infty}{(q;q)_n}
(-c)^n q^{n(n+k)},
\]
which is exactly the left-hand side of \eqref{eq:prop}.
The special case $c=-1$ can also be obtained by reversing the order of
summation.
\end{proof}

From \eqref{eq:resolution} and \eqref{eq:lj} it follows that
\begin{align*}
\bigl< E\bigl(\{q^n\}\bigr)v,w \bigr>&= \frac{-1}{2\pi
i}\oint_{(q^n)}
\bigl< (J^\ast-s)^{-1}v,w \bigr>\,ds \\
&=\frac{-1}{2\pi i}\sum_{l\leq
j}(v_l\overline{w}_j+v_j\overline{w}_l)
(1-\tfrac{1}{2}\delta_{j,l})\oint_{(q^n)}
\frac{\psi_j(s)\Psi_l(s)}{[\psi(s),\Psi(s)]}\,ds.
\end{align*}
The integral on the right-hand side is given by
\[
\frac{-1}{2\pi i}\oint_{(q^n)}
\frac{\psi_j(s)\Psi_l(s)}{[\psi(s),\Psi(s)]}\,ds=
\psi_j(q^n)\Psi_l(q^n)\,\underset{z=q^n}{\res}\frac{1}{[\psi(z),\Psi(z)]}
\]
and by Lemma \ref{lem:proportional} (with $c=-1$), we have
$\psi_k(q^n)=(-1)^nt^nq^{n^2}\Psi_k(q^n)$. Combining this with the
fact that
\begin{align*}
\underset{z=q^n}{\res}\frac{1}{[\psi(z),\Psi(z)]}
=\frac{(-1)^{n+1}t^nq^{n(n+1)}}{(q;q)_n}\frac{1}{(-tq,-1/t,q;q)_\infty},
\end{align*}
we end up with
\[
\bigl< E\bigl(\{q^n\}\bigr)v,w \bigr>=
\frac{q^n}{(q;q)_n}\frac{\bigl<v,\psi(q^n)\bigr>\,
\bigl<\psi(q^n),w\bigr>}{(-tq,-1/t,q;q)_\infty}.
\]
In particular, it follows that
\begin{equation}
\label{eq:q^n}
\Vert \psi(q^n) \Vert^2=\frac{(q;q)_n}{q^n}(-tq,-1/t,q;q)_\infty
\quad\mbox{and}\quad
\bigl< \psi(q^n),\psi(q^m) \bigr>=\Vert \psi(q^n)\Vert^2\delta_{m,n}
\end{equation}
if we set $v=w=\psi(q^n)$, respectively $v=w=\psi(q^m)$.

In a similar way as above, one can show that
\[
\bigl< E\bigl(\{-q^r/t\}\bigr)v,w \bigr>=
\frac{q^r}{(-q/t;q)_r}\frac{\bigl<v,\psi(-q^r/t)\bigr>\,
\bigl<\psi(-q^r/t),w\bigr>}{(-t,q,q;q)_\infty}.
\]
For by Lemma \ref{lem:proportional}, we have
$\psi_k(-q^r/t)=(-1)^rq^{r^2}t^{-r}\Psi_k(-q^r/t)$ and
\[
\underset{z=-q^r/t}{\res}\frac{1}{[\psi(z),\Psi(z)]}
=\frac{(-1)^{r+1}q^{r(r+1)}}{t^r(-q/t;q)_r}\frac{1}{(-t,q,q;q)_\infty}.
\]
It thus follows that
\begin{equation}
\label{eq:-q^r/t}
\bigl<\psi(-q^r/t),\psi(-q^s/t)\bigr>=
\frac{(-q/t;q)_r}{q^r}(-t,q,q;q)_\infty\delta_{r,s}.
\end{equation}
Moreover, we clearly have
\begin{equation}
\label{eq:mix}
\bigl<\psi(q^n),\psi(-q^r/t)\bigr>=0.
\end{equation}

Recall now that the Stieltjes--Wigert polynomials are given by
\[
S_n(x;q)=\frac{1}{(q;q)_n}\phi_{q^n}(x)
\]
and consider also the functions $M_r^{(t)}(x;q)$ defined by
\[
M_r^{(t)}(x;q)=\frac{1}{(q;q)_\infty}\phi_{-q^r/t}(x), \quad
r\in\Z.
\]
These functions are closely related to the second $q$-Bessel function
\cite[Exerc.~1.24]{GaspR} defined by
\[
J_\nu^{(2)}(z;q)=\frac{(z/2)^\nu}{(q;q)_\infty}
\rps{1}{1}{-z^2/4}{0}{q^{\nu+1}}
=\frac{(q^{\nu+1};q)_\infty}{(q;q)_\infty}(z/2)^\nu
\rps{0}{1}{-}{q^{\nu+1}}{-\frac{z^2q^{\nu+1}}{4}}.
\]
Indeed, we have
$t^{\frac{k+r}{2}}M_r^{(t)}(tq^k;q)=
q^{r(r+k)/2}J_{k+r}^{(2)}(2\sqrt{t}q^{-r/2};q)$.

It follows immediately from Proposition \ref{prop:eigenfunctions} that
\[
L S_n(\,\cdot\,;q)=q^nS_n(\,\cdot\,;q) \quad\mbox{for } n\in\Z_+
\]
and
\[
LM_r^{(t)}(\,\cdot\,;q)=-\frac{q^r}{t} M_r^{(t)}(\,\cdot\,;q)
\quad \mbox{for } r\in\Z.
\]
Furthermore, since the spectral decomposition is unique, these
eigenfunctions form an orthogonal basis for $L^2(\mu_t)$. We put
together the results from \eqref{eq:q^n}, \eqref{eq:-q^r/t} and
\eqref{eq:mix} in the following theorem which is a formal limit
transition of \cite[Thm.~4.1]{CiccKK}.
\begin{thm}
\label{OR}
The Stieltjes--Wigert polynomials $S_n(x;q)$, respectively the
$q$-Bessel functions $M_r^{(t)}(x;q)$, are orthogonal in $L^2(\mu_t)$.
The orthogonality relations are given by
\begin{equation}
\label{eq:OR-Sn}
\frac{1}{(-tq,-1/t,q;q)_\infty}
\sum_{k=-\infty}^\infty t^kq^{\binom{k+1}{2}}
S_n(tq^k;q)S_m(tq^k;q)=\frac{\delta_{m,n}}{q^n(q;q)_n}
\end{equation}
and
\begin{equation}
\label{eq:OR-Mr}
\frac{1}{(-t;q)_\infty}\sum_{k=-\infty}^\infty t^kq^{\binom{k+1}{2}}
M_r^{(t)}(tq^k;q)M_s^{(t)}(tq^k;q)=
\frac{(-q/t;q)_r}{q^r}\delta_{r,s}.
\end{equation}
Moreover, $S_n(x;q)$ and $M_r^{(t)}(x;q)$ are mutually orthogonal in
$L^2(\mu_t)$, that is,
\begin{equation}
\label{eq:OR-MS}
\sum_{k=-\infty}^\infty t^kq^{\binom{k+1}{2}}
S_n(tq^k;q)M_r^{(t)}(tq^k;q)=0 \quad \mbox{for all $n$, $r$}
\end{equation}
and $\bigl\{S_n(x;q)\bigr\}_{n\in\Z_+} \cup
\bigl\{M_r^{(t)}(x;q)\bigr\}_{r\in\Z}$ form an orthogonal basis for
$L^2(\mu_t)$.
\end{thm}

\begin{remark}
The orthogonality relation \eqref{eq:OR-Sn} is due to Chihara
\cite{Chihara}, whereas \eqref{eq:OR-Mr} is the Hansen--Lommel
orthogonality relation for the second $q$-Bessel function, see
\cite[Thm.~3.1]{KoelJMAA}. The above theorem contradicts
\cite[Thm.~3.3]{KoelJMAA}, and the flaw in the proof of
\cite[Thm.~3.3]{KoelJMAA} is contained in \cite[Lemma 3.4]{KoelJMAA},
where the unbounded operator $S$ as constructed there is not symmetric
as claimed.

The statement in \eqref{eq:OR-MS} can also be proved directly in the
following way. Use \cite[(0.6.8/9)]{K&S} to write $M_r^{(t)}(x;q)$ as
\[
M_r^{(t)}(x;q)=\frac{(xq^{r+1}/t;q)_\infty}{(q;q)_\infty}
\rps{0}{1}{-}{xq^{r+1}/t}{-xq},
\]
so that
\begin{multline*}
\sum_{k=-\infty}^\infty t^kq^{\binom{k+1}{2}}
S_n(tq^k;q)M_r^{(t)}(tq^k;q)=\\
\sum_{k=-\infty}^\infty t^kq^{\binom{k+1}{2}}
\frac{1}{(q;q)_n}\rps{1}{1}{q^{-n}}{0}{-tq^{k+n+1}}
\frac{(q^{k+r+1};q)_\infty}{(q;q)_\infty}\rps{0}{1}{-}{q^{k+r+1}}{-tq^{k+1}}.
\end{multline*}
Because of absolute convergence we can interchange the order of
summation to get
\[
\frac{1}{(q;q)_n}\sum_{m=0}^n\frac{(q^{-n};q)_m}{(q;q)_m}
t^mq^{\binom{m}{2}+m(n+1)}
\frac{1}{(q;q)_\infty}\sum_{l=0}^\infty\frac{(-1)^lq^{l^2}t^l}{(q;q)_l}
\sum_{k=-\infty}^\infty(q^{k+r+l+1};q)_\infty t^kq^{\binom{k}{2}+k(m+l+1)}.
\]
The inner sum (over $k$) reduces to
\begin{align*}
\sum_{k=-r-l}^\infty (q^{k+r+l+1};q)_\infty
t^kq^{\binom{k}{2}+k(m+l+1)}&=
\frac{(q;q)_\infty q^{\binom{r+l}{2}}}{t^{r+l}q^{(r+l)(m+l)}}
\sum_{k=0}^\infty\frac{t^kq^{\binom{k}{2}+k(m+1-r)}}{(q;q)_k} \\
&=\frac{(-tq^{m+1-r},q;q)_\infty q^{\binom{r}{2}+\binom{l}{2}}}
{q^{l^2+m(r+l)}t^{r+l}}
\end{align*}
and the sum over $l$ then becomes
\[
\sum_{l=0}^\infty\frac{(-1)^lq^{\binom{l}{2}-lm}}{(q;q)_l}=
(q^{-m};q)_\infty.
\]
Since $(q^{-m};q)_\infty=0$ for $m\geq 0$, the relation
\eqref{eq:OR-MS} is established.
\end{remark}

\begin{remark}
Using the explicit expression for $M_r^{(t)}(x;q)$ and Lemma
\ref{lem:proportional}, we see that $|M_r^{(t)}(tq^k;q)|$ is bounded by
some constant, say $M(r,t)$, for all $k\in\Z$ provided $t<q^r$. By the
construction of Berg \cite{Berg-geo} it thus follows from Theorem
\ref{OR} that the measure
\[
\nu_{s,t}=\frac{1}{(-tq,-1/t,q;q)_\infty}\sum_{k=-\infty}^\infty
t^kq^{\binom{k+1}{2}} \Bigl( 1+\frac{s}{M(r,t)}{M_r^{(t)}(tq^k;q)}
\Bigr)\delta_{tq^k}
\]
is a solution to the Stieltjes--Wigert moment problem for all $|s|\leq
1$ and $t<q^r$.
\end{remark}

\section{Spectral analysis for the continuous case}
\label{continuous}
We now work on the Hilbert space $L^2(\mu)$, where $\mu$ is the
absolutely continuous measure from Corollary
\ref{cor:symmetricmeasure} (ii). The density of $\mu$, which will be
denoted $w$, thus satisfies the functional equation
\begin{equation}
\label{eq:func-eq}
w(xq)=xw(x), \quad x>0.
\end{equation}
We remind the reader that a function $g$ is called $q$-periodic if
$g(xq)=g(x)$ for all $x>0$.

\subsection{Direct integral decomposition}

Consider the Hilbert space $\ell^2(\Z)$ equipped with its standard
orthonormal basis $\{ e_k\}_{k\in\Z}$. For a compactly supported
measurable function $f$ on $(0,\infty)$ we define
\begin{align}
\label{I}
\notag
(q,1]\ni t\mapsto (If)(t)&=\sum_{k=-\infty}^\infty f(tq^k) q^{k/2}
\sqrt{w(tq^k)}\,e_k  \\
&=\sqrt{w(t)}\sum_{k=-\infty}^\infty f(tq^k)t^{k/2}q^{k(k+1)/4}
\,e_k  \in \ell^2(\Z).
\end{align}
Clearly, $(I(gf))(t)=g(t)(If)(t)$ whenever $g$ is a $q$-periodic
function.

\begin{prop}
\label{prop:L2muisdirectintegral}
The operator $I$ defined in \eqref{I} extends to a unitary isomorphism
\[
I\colon L^2(\mu) \rightarrow \int_\Om^{\oplus}\ell^2(\Z)\,dt
\]
with $\Om=(q,1]\cap \text{\rm supp}(\mu)$.
\end{prop}

\begin{remark}
\label{rmk:propL2muisdirectintegral}
The direct integral Hilbert space $\int_\Om^{\oplus}\ell^2(\Z)\,dt$
consists of all measurable functions $f: \Om\ra\ell^2(\Z)$ with
$\int_\Om\|f(t)\|^2_{\ell^2(\Z)}\,dt<\infty$. The term measurable
means that $t\mapsto\langle f(t), e_k \rangle_{\ell^2(\Z)}$ is
measurable for all $k\in\Z$. In particular, the constant vector fields
$t\mapsto e_j$ are measurable. The inner product on
$\int_\Om^{\oplus}\ell^2(\Z)\,dt$ is given by
\[
\langle f, g \rangle_{\int_\Om^{\oplus}\ell^2(\Z)\,dt}=
\int_\Om\langle f(t), g(t) \rangle_{\ell^2(\Z)}\,dt
\]
and we have $\int_\Om^\oplus\ell^2(\Z)\,dt \cong L^2(\Om)\otimes
\ell^2(\Z)$ as Hilbert spaces. The space of all $t\mapsto g(t)e_j$,
$g$ bounded measurable function on $\Om$, is therefore dense in
$\int_\Om^{\oplus}\ell^2(\Z)\,dt$.
Notice that $t\mapsto (If)(t)$ as defined in \eqref{I} is measurable.
See e.g. \cite[Part II, Ch. 1]{Dixm} for more information.
\end{remark}

\begin{proof}
For $f,g$ compactly supported functions in $L^2(\mu)$, we have
\begin{equation*}
\begin{split}
\langle If, Ig \rangle_{\int_\Om^\oplus
\ell^2(\Z) dt} &= \int_\Om \langle (If)(t),
(Ig)(t) \rangle_{\ell^2(\Z)} dt
= \int_\Om \sum_{k=-\infty}^\infty f(tq^k) \overline{g(tq^k)}
q^k w(tq^k) \, dt \\
&= \sum_{k=-\infty}^\infty \int_q^1 f(tq^k) \overline{g(tq^k)}
q^k w(tq^k)\, dt = \sum_{k=-\infty}^\infty
\int_{q^{k+1}}^{q^k} f(x) \overline{g(x)} w(x)\, dx \\ &=
\int_0^\infty f(x)\overline{g(x)} w(x)\, dx =
\langle f, g\rangle_{L^2(\mu)},
\end{split}
\end{equation*}
where interchanging summation and integration is allowed since $f,g$
being compactly supported implies that the sum is finite. Moreover, we
can switch from $\int_\Om$ to $\int_q^1$ since $w$ satisfies the
functional equation \eqref{eq:func-eq}.

Recalling that the compactly supported measurable functions are dense
in $L^2(\mu)$, the operator $I$ from \eqref{I} extends to an isometry
$I\colon L^2(\mu)\to\int_\Om^\oplus \ell^2(\Z)dt$. Since the image of
$I$ contains any element of the form $t\mapsto h(t)e_k$, $h$ bounded
measurable function on $\Om$, and these elements are dense in
$\int_\Om^\oplus \ell^2(\Z)dt$, we conclude that $I\colon
L^2(\mu)\to\int_\Om^\oplus \ell^2(\Z)dt$ is surjective and thus
unitary.
\end{proof}

The adjoint of the unitary operator $I$ is given explicitly by
\begin{equation}
\label{eq:adjointI}
I^\ast \Bigl( t\mapsto \sum_{k=-\infty}^\infty h_k(t)\,e_k\Bigr)(x)
=\sum_{k=-\infty}^\infty \chi_{(q^{k+1},q^k]}(x)
\frac{h_k(xq^{-k})}{q^{k/2} \sqrt{w(x)}},
\end{equation}
where $\chi_A$ denotes the indicator function of the set $A$. The
right-hand side of \eqref{eq:adjointI} only makes sense when $w(x)>0$,
but there is no need to specify the value of a function in $L^2(\mu)$
at points where $w(x)=0$. Formally calculating $I\phi_z$, with
$\phi_z$ the eigenfunction of $L$ from Proposition
\ref{prop:eigenfunctions}, gives
\begin{equation*}
(I\phi_z)(t) = \sqrt{w(t)}
\sum_{k=-\infty}^\infty \phi_z(tq^k)t^{k/2}q^{k(k+1)/4}\,e_k
= \sqrt{w(t)}\,\psi(z;t),
\end{equation*}
with $\psi(z;t)$ the formal, i.e. in general not contained in
$\ell^2(\Z)$, eigenvectors of $J_t$ as in Section \ref{WandG}.
Conversely, by \eqref{eq:adjointI} we have for any function $f$ on
$\Om$ that
\begin{equation*}
I^\ast\Bigl( t\mapsto f(t) \sum_{k=-\infty}^\infty \phi_z(tq^k)
t^{k/2}q^{k(k+1)/4} e_k\Bigr) =
\text{Per}({f}/{\sqrt{w}})\,\phi_z,
\end{equation*}
where $\text{Per}$ maps a function on $\Om$ to a $q$-periodic function
on $\text{supp}(\mu)$ such that they are equal on $\Om$, explicitly
\begin{equation}\label{eq:defPer}
\text{Per}(f)(x)=\sum_{k=-\infty}^\infty \chi_{(q^{k+1}, q^k]}(x) f(xq^{-k}).
\end{equation}

Recall from Section \ref{ell2} the unbounded symmetric operator
$(J_t,\mathcal{D})$ on $\ell^2(\Z)$ defined by
\[
J_t e_k = a_k(t) e_{k+1} + b_k(t) e_k + a_{k-1}(t) e_{k-1},
\quad k\in\Z
\]
with
\[
a_k(t)=\frac{1}{\sqrt{tq^{k+1}}}
\quad\mbox{and}\quad
b_k(t)=-\frac{1}{tq^k}.
\]
Note that $a_k$ and $b_k$ are bounded continuous functions of $t\in
(q,1]$ for fixed $k\in\Z$.  It follows from Theorem
\ref{thm:self-adjoint} that $(J_t,\mathcal{D})$ is essentially
self-adjoint, and we denote by $(J_t^\ast, \text{dom}(J_t^\ast))$ its
unique self-adjoint extension.

Let $L^2(\Om)\otimes {\mathcal D}$ be the (algebraic) tensor product
of the space $L^2(\Om)$ and the space ${\mathcal D}$ of finite linear
combinations of the basis vectors. By Remark
\ref{rmk:propL2muisdirectintegral} this tensor product is dense in
$\int_\Om^\oplus \ell^2(\Z)dt$ since it contains $B(\Om)\otimes
{\mathcal D}$, with $B(\Om)$ the space of bounded measurable functions
on $\Om$. Observe that for $h\otimes v\in L^2(\Om)\otimes {\mathcal
  D}$, the field $t\mapsto h(t) J_tv$ is measurable because the inner
product
\[
t\mapsto \langle h(t) J_tv, e_k\rangle= h(t) \langle v, J_te_k\rangle
= h(t)\bigl(a_k(t) \langle v,e_{k+1}\rangle + b_k(t) \langle
v,e_k\rangle + a_{k-1}(t) \langle v,e_{k-1}\rangle\bigr)
\]
is measurable for any $k\in\Z$. Moreover, this inner product is only
non-zero for finitely many values of $k$, so the vector field
$t\mapsto h(t) J_tv$ is an element of $\int_\Om^\oplus \ell^2(\Z)dt$.
We now define $\int_\Om^\oplus J_t dt$ as the operator with domain
$L^2(\Om)\otimes {\mathcal D}$ mapping the element $h\otimes v$
considered as the vector field $t\mapsto h(t)v$ to $t\mapsto h(t)
J_tv$. Note that $h\otimes v$ is identified with $f\otimes v$ whenever
$f=h$ a.e. in $\Om$.

\begin{prop}
\label{prop:Lindirectintegral}
Consider $L$ as an unbounded operator with domain the compactly
supported functions in $L^2(\mu)$.
Then $I$ intertwines $L$ with $J=\int_\Om^\oplus J_t\,dt$.
\end{prop}

\begin{proof}
For $f$ compactly supported, take $N,M\in\Z$ such that
$\text{supp}(f)\subset (q^{N+1},q^M]$ and identify
\[
(If)(t) = \sum_{k=N}^M f(tq^k) q^{k/2} \sqrt{w(tq^k)}\, e_k
\]
with $\sum_{k=N}^M h_k\otimes e_k \in L^2(\Om)\otimes {\mathcal D}$,
where $h_k(t)=f(tq^k)q^{k/2} \sqrt{w(tq^k)}$. Since
\[
\int_\Om |h_k(t)|^2\, dt = \int_{q^{k+1}}^{q^k}
|f(x)|^2w(x)dx<\infty,
\]
we have indeed $h_k\in L^2(\Om)$.  So $I$ maps the domain of $L$ into
$L^2(\Om)\otimes {\mathcal D}$. Conversely, $I^\ast$ of an element
$h\otimes e_k\in L^2(\Om)\otimes {\mathcal D}$ gives by
\eqref{eq:adjointI} a compactly supported function on $(0,\infty)$ and
\[
\int_0^\infty |I^\ast(h\otimes e_k)(x)|^2w(x)dx =
\int_\Om |h(t)|^2dt<\infty.
\]
The intertwining property is a straightforward calculation. For
$f\in\text{dom}(L)$ and fixed $t\in\Om$, we have
\begin{equation*}
\begin{split}
I(Lf)(t) &=
\sqrt{w(t)}\sum_{k=-\infty}^\infty
\Bigl(f(tq^{k+1})-\frac{1}{tq^k}f(tq^k)+\frac{1}{tq^k}f(tq^{k-1})\Bigr)
t^{k/2}q^{k(k+1)/4} e_k \\
&= \sqrt{w(t)}\sum_{k=-\infty}^\infty f(tq^k)
\Bigl(\frac{1}{\sqrt{tq^k}}\,e_{k-1}-\frac{1}{tq^k}\,e_k +
\frac{1}{\sqrt{tq^{k+1}}}\,e_{k+1}\Bigr)t^{k/2}q^{k(k+1)/4}
= J_t (If)(t).
\end{split}
\end{equation*}
Note that the infinite sums only contain a finite number of non-zero
terms, so that all rearrangements are valid.
\end{proof}

Since the operator $L$ from Proposition \ref{prop:Lindirectintegral}
is symmetric and commutes with complex conjugation, it has a
self-adjoint extension. We aim at finding its adjoint for which we
want to give a direct integral representation. Because of Proposition
\ref{prop:Lindirectintegral} and the fact that each
$(J_t^\ast,\text{dom}(J_t^\ast))$ is self-adjoint we consider the
operator $J^\ast = \int_\Om^\oplus J_t^\ast dt$. The next paragraph
justifies this notation.

According to \cite[Def. p.~283]{ReedS4} we need to check that the field
of operators $t\mapsto (J_t^\ast+i)^{-1}$ is measurable, i.e. that
$t\mapsto \langle (J_t^\ast+i)^{-1}e_k,e_l\rangle_{\ell^2(\Z)}$ is
measurable for all $k,l\in\Z$. By the functional calculus for
$J_t^\ast$ established in Section \ref{discrete}, we have
\[
\langle (J_t^\ast+i)^{-1}e_k,e_l\rangle_{\ell^2(\Z)}=
\int_\R \frac{1}{\la+i} dE^t_{e_k,e_l}(\la),
\]
where the right-hand side can be written as
\[
\sum_{n=0}^\infty \frac{1}{q^n+i}\frac{\langle e_k,\psi(q^n;t)\rangle
\langle\psi(q^n;t),e_l\rangle}{\|\psi(q^n;t)\|^2}
+\sum_{r=-\infty}^\infty \frac{1}{i-q^r/t}
\frac{\langle e_k,\psi(-q^r/t;t)\rangle\langle \psi(-q^r/t;t),e_l\rangle}
{\|\phi(-q^r/t;t)\|^2}.
\]
The desired measurability hence follows.
Now define
\begin{equation*}
\begin{split}
\text{dom}(J^\ast) 
&=\Bigl\{ t\mapsto u(t) \in \int_\Om^\oplus \ell^2(\Z)dt
\,\Big|\, u(t) \in \text{dom}(J_t^\ast)\text{ a.e.},
\int_\Om \| J_t^\ast u(t)\|^2dt <\infty \Bigr\}, \\
&J^\ast = \int_\Om^\oplus J_t^\ast dt \colon
\text{dom}(J^\ast) 
\ni \bigl(t\mapsto u(t)\bigl) \longmapsto
\bigl(t\mapsto J_t^\ast u(t)\bigl).
\end{split}
\end{equation*}
By \cite[Thm. XIII.85, p.~284]{ReedS4} the operator
$J^\ast=\int_\Om^\oplus J_t^\ast dt$ is the adjoint of $J$ and
$J^\ast$ is self-adjoint. Moreover, the functional calculus is given
by
\begin{equation}\label{eq:FCfordirectintegral}
f(J^\ast) = f\Bigl( \int_\Om^\oplus J_t^\ast dt\Bigr) =
\int_\Om^\oplus f(J_t^\ast) dt
\end{equation}
for any bounded measurable function $f$ on $\R$.

\begin{prop}\label{prop:adjointLasdirectintegral}
  The adjoint operator $(L^\ast,\text{\rm dom}(L^*))$ is intertwined
  with $(J^\ast, \text{\rm dom}(J^\ast))$ by the unitary isomorphism
  $I$.
\end{prop}

As an immediate consequence, we have
\begin{cor}\label{cor:prop:adjointLasdirectintegral}
  $(L^\ast,\text{\rm dom}(L^\ast))$ is the unique self-adjoint
  extension of $(L,\text{\rm dom}(L))$, and for any bounded Borel
  function $f$ on $\R$ the functional calculus is given by
\begin{equation*}
f(L^\ast) = I^\ast \int_\Om^\oplus f(J_t^\ast) dt\, I.
\end{equation*}
\end{cor}

\begin{proof}
[Proof of Proposition \ref{prop:adjointLasdirectintegral}]
The domain of $L^\ast$ consists of all functions $g\in L^2(\mu)$ such
that
\begin{equation*}
\begin{split}
f \mapsto \langle Lf,g\rangle_{L^2(\mu)}
= \int_\Om \langle I(Lf)(t), (Ig)(t)\rangle_{\ell^2(\Z)} dt
= \int_\Om \langle J_t(If)(t), (Ig)(t)\rangle_{\ell^2(\Z)} dt
\end{split}
\end{equation*}
defines a continuous linear functional on $\text{dom}(L)$. We have
used Proposition \ref{prop:Lindirectintegral} to replace $I(Lf)$ with
$J_t(If)$ in the inner product on the right-hand side. So for
$g\in\text{dom}(L^\ast)$ there exists a constant $C=C(g)>0$ such that
\begin{equation}\label{eq:pfprop:adjointLasdirectintegral1}
|\langle Lf,g\rangle_{L^2(\mu)}|
=\Bigl|\int_\Om \langle J_t(If)(t), (Ig)(t)\rangle_{\ell^2(\Z)} dt\Bigr|
\leq C \|f\|_{L^2(\mu)} = C
\Bigl( \int_\Om \|If(t)\|_{\ell^2(\Z)}^2 dt\Bigr)^{1/2}
\end{equation}
for all $f\in\text{dom}(L)$. Since $f$ is compactly
supported, the inner product $\langle J_t(If)(t),
(Ig)(t)\rangle_{\ell^2(\Z)}$ is a finite sum and hence equal to
$\langle If(t), J_t^\ast(Ig)(t)\rangle_{\ell^2(\Z)}$. Therefore,
\eqref{eq:pfprop:adjointLasdirectintegral1} can be rewritten as
\begin{equation}\label{eq:pfprop:adjointLasdirectintegral2}
\Bigl| \int_\Om \langle If(t), J_t^\ast (Ig)(t)\rangle_{\ell^2(\Z)}
dt \Bigr| \leq C \Bigl(\int_\Om \| If(t)\|_{\ell^2(\Z)}^2 dt\Bigr)^{1/2}.
\end{equation}
Now we can show that the vector field $t\mapsto Ig(t)$ belongs to
$\text{dom}(\int_\Om^\oplus J_t^\ast dt)$ whenever
$g\in\text{dom}(L^\ast)$. First, by taking $f= I^\ast (1\otimes
e_k)\in \text{dom}(L)$ we see that
\begin{equation*}
t\mapsto \langle e_k, J_t^\ast(Ig)(t)\rangle_{\ell^2(\Z)}
= \langle If(t), J_t^\ast(Ig)(t)\rangle_{\ell^2(\Z)}
= \langle J_t(If)(t), Ig(t)\rangle_{\ell^2(\Z)}
\end{equation*}
is measurable and square integrable on $\Om$ for any $k\in\Z$, since
\[
\langle J_t(If)(t), Ig(t)\rangle_{\ell^2(\Z)}=
a_k(t) \langle e_{k+1}, Ig(t)\rangle_{\ell^2(\Z)}
+ b_k(t) \langle e_k, Ig(t)\rangle_{\ell^2(\Z)}
+ a_{k-1}(t) \langle e_{k-1}, Ig(t)\rangle_{\ell^2(\Z)}.
\]
Then apply \eqref{eq:pfprop:adjointLasdirectintegral2} with $If(t) =
\sum_{k=-N}^N \langle J_t^\ast (Ig)(t),e_k\rangle_{\ell^2(\Z)} e_k$
to get
\begin{equation*}
\sum_{k=-N}^N \int_\Om |\langle J_t^\ast (Ig)(t),e_k\rangle_{\ell^2(\Z)}|^2
dt \leq  C \Bigl(\int_\Om \sum_{k=-N}^N
|\langle J_t^\ast (Ig)(t),e_k\rangle_{\ell^2(\Z)} |^2 dt\Bigr)^{1/2}
\end{equation*}
or
\begin{equation*}
\sum_{k=-N}^N \int_\Om |\langle J_t^\ast (Ig)(t),e_k\rangle_{\ell^2(\Z)}|^2
dt \leq  C^2.
\end{equation*}
Since $C$ is independent of $N$, this is also valid for
$N\to\infty$. In particular, it follows that
\[
\sum_{k=-\infty}^\infty
|\langle J_t^\ast (Ig)(t),e_k\rangle_{\ell^2(\Z)}|^2 <\infty
\quad \mbox{a.e.}
\]
so that $t\mapsto J_t^\ast (Ig)(t)$ is a measurable square integrable
vector field for which $Ig(t) \in \text{dom}(J_t^\ast)$ a.e.  This
proves that $I \text{dom}(L^\ast) \subset \text{\rm
  dom}(\int_\Om^\oplus J_t^\ast dt)$ and $IL^\ast$ is the restriction
of $J^\ast I=\int_\Om^\oplus J_t^\ast dt\, I$.

For the converse inclusion take $g\in I^\ast\text{\rm dom}(J^\ast)$
and observe that for any $f\in\text{dom}(L)$,
\begin{equation*}
\begin{split}
|\langle Lf,g\rangle_{L^2(\mu)}| =&
\Bigl|\int_\Om \langle J_t(If)(t),Ig(t)\rangle_{\ell^2(\Z)} dt\Bigl|
= \Bigl|\int_\Om \langle If(t),J_t^\ast(Ig)(t)\rangle_{\ell^2(\Z)} dt\Bigl| \\
\leq & \Bigl(\int_\Om \| If(t)\|^2_{\ell^2(\Z)} dt\Bigr)^{1/2}
\Bigl(\int_\Om \| J_t^\ast(Ig)(t)\|^2_{\ell^2(\Z)} dt \Bigr)^{1/2}
= C \| f\|_{L^2(\mu)}.
\end{split}
\end{equation*}
In other words, $f\mapsto \langle Lf,g\rangle_{L^2(\mu)}$ defines a
continuous linear functional on $\text{dom}(L)$ and it follows that
$I^\ast\text{dom}(J^\ast)\subset\text{dom}(L^\ast)$.
\end{proof}

\subsection{Spectral decomposition for $L^\ast$}

We start this section by presenting the spectrum of $L^\ast$.

\begin{thm} \label{thm:spectraldecompL}
The spectrum of the self-adjoint operator $(L^\ast,\text{\rm
  dom}(L^\ast))$ consists of point spectrum $q^{\Z_+}$, each
point having infinite multiplicity, and continuous spectrum
$\overline{\cup_{l\in\Z} \tilde \Om_l}$, where $\tilde\Om_l = \{
-q^l/t \mid t\in\Om\}$. In particular, we have $\si(L^\ast) =
(-\infty,0]\cup q^{\Z_+}$ when $\Om=(q,1]$.
\end{thm}
\begin{proof}
The theorem follows from \cite[Thm.~XIII.85]{ReedS4} and Proposition
\ref{prop:adjointLasdirectintegral}.  We only need to consider the
point $0$ which is in the closure of $q^{\Z_+}$ and in the
closure of $\cup_{l\in\Z} \tilde\Om_l$.  Since $(L^\ast,\text{\rm
  dom}(L^\ast))$ is self-adjoint, $0$ is either in the point spectrum
or in the continuous spectrum. In case $0$ is in the point spectrum,
it is also contained in the point spectrum of $(J^\ast,
\text{dom}(J^\ast))$, so there exists a non-trivial function $t\mapsto
v(t)$ such that $J_t^\ast v(t)=0$ a.e. on $\Om$. By Theorem
\ref{thm:spectrum}, however, the point $0$ is not contained in the
point spectrum of $(J_t^\ast, \text{dom}(J_t^\ast))$ for any $t\in
\Om$, so $v(t)=0$ a.e. and $0$ belongs to the continuous spectrum.
\end{proof}

In order to make Theorem \ref{thm:spectraldecompL} more explicit we
establish the corresponding spectral decomposition. Following the
ideas of the proof of \cite[Thm.~XIII.86]{ReedS4} we define
\begin{equation}
\label{eq:defHn}
{\mathcal H}^+_n = \Bigl\{ v\in \int_\Om^\oplus \ell^2(\Z)dt
\,\Big|\, v(t)=f(t) \frac{\psi(q^n;t)}{N_{q^n}(t)}
\text{ for some } f\in L^2(\Om)\Bigr\},
\quad n\in\Z_+,
\end{equation}
and
\begin{equation}
\label{eq:defHr}
{\mathcal H}^-_r = \Bigl\{ v\in \int_\Om^\oplus \ell^2(\Z)dt \,\Big|\,
v(t) = f(t) \frac{\psi(-q^r/t;t)}{N_{-q^r/t}(t)}
\text{ for some } f\in L^2(\Om)\Bigr\},
\quad r\in\Z,
\end{equation}
using the notation $N_{\xi}(t)=\|\psi(\xi;t)\|_{\ell^2(\Z)}$ for $\xi$
in the point spectrum of $J_t^\ast$. Then ${\mathcal H}_n^+$,
${\mathcal H}^-_r$ are mutually orthogonal closed subspaces of
$\int_\Om^\oplus \ell^2(\Z)dt$ and, moreover,
\[
\int_\Om^\oplus\ell^2(\Z)dt = {\mathcal H}^+\oplus {\mathcal H}^-,
\quad \mbox{with }
{\mathcal H}^+ = \bigoplus_{n=0}^\infty {\mathcal H}_n^+
\mbox{ and }
{\mathcal H}^- = \bigoplus_{r=-\infty}^\infty {\mathcal H}^-_r.
\]
Note that the subspaces ${\mathcal H}^\pm_l$ are contained in
$\text{dom}(J^\ast)$ and $J^\ast$ preserves each of them. By $U^\pm_l
\colon{\mathcal H}^\pm_l \to L^2(\Om)$ we denote the unitary operator
defined by $U^\pm_l v = f$ for $v\in {\mathcal H}^\pm_l$ of the form
as in \eqref{eq:defHn} or \eqref{eq:defHr}. It follows that $U^\pm_l$
intertwines $J^\ast$ with multiplication by $\la^\pm_l$ on $L^2(\Om)$,
where $\la^+_l(t)=q^l$ and $\la^-_l(t)=-q^l/t$. We put $J_l^\pm =
U_l^\pm J^\ast (U_l^\pm)^\ast$ so that $J^\pm_l f = \la^\pm_l f$ for
all $f\in L^2(\Om)$. In particular, it follows that $\text{ker}(J^\ast
-q^l)= {\mathcal H}^+_l$ so that $q^{\Z_+}$ is contained in the
point spectrum of $J^\ast$, and each point of this form has infinite
multiplicity.

For the case of negative eigenvalues we define $\tilde\Om_l = \{
-q^l/t \mid t\in\Om\} \subseteq (-q^{l-1},-q^l]$ for $l\in\Z$. Then
$V_l \colon L^2(\Om) \to L^2(\tilde\Om_l)$ given by
\[
(V_lf)(\la) = \frac{q^{l/2}}{|\la|} f(-q^l/\la),
\quad \la\in\tilde{\Om}_l
\]
is a unitary operator and its adjoint $V_l^\ast$ is almost given by
the same formula,
\[
(V_l^\ast g)(t)=\frac{q^{l/2}}{t}g(-q^l/t),
\quad t\in\Om.
\]
By a straightforward calculation we see that
\begin{equation}\label{eq:multoperrealJ-l}
(V_l^{} J^-_l V_l^\ast g)(\la) = \la \, g(\la),
\quad \la\in\tilde{\Om}_l
\end{equation}
for any $g\in L^2(\tilde{\Om}_l)$. It thus follows that $\tilde \Om
=\cup_{l\in\Z} \tilde\Om_l \subseteq (-\infty, 0]$ is contained in the
continuous spectrum of $J^\ast$, and this part of the spectrum is
simple. Using the notation $E(T|A)$ for the spectral projection
corresponding to the Borel set $A\subset\R$ for a (possibly unbounded)
self-adjoint operator $T$, we see that $E(V^{}_l J^-_l V_l^\ast |
A)$ is just multiplication by the characteristic function $\chi_{A\cap
  \tilde\Om_l}$. Tracing the steps back it follows that
\[
E(J^\ast\vert_{{\mathcal H}^-_l}| A) v(t) = \chi_{A\cap
  \tilde\Om_l}(-q^l/t) v(t),
\]
with the notation as in \eqref{eq:defHn} and \eqref{eq:defHr}. By
considering $J^\ast$ restricted to ${\mathcal H}^-$, we see that
$\si(J^\ast\vert_{{\mathcal H}^-}) = \overline{\cup_{l\in\Z} \tilde
  \Om_l}$.

To obtain the spectral decomposition $E$ of $(L^\ast,
\text{dom}(L^\ast))$ we use Proposition
\ref{prop:adjointLasdirectintegral} and Theorem
\ref{thm:spectraldecompL}. The idea is to get the results from the
spectral decomposition for $J^\ast$ using the unitary isomorphism $I$.
First we consider the spectral decomposition corresponding to the
point spectrum $\sigma_p(L^\ast)$. It follows that $L^\ast$ preserves
$I^\ast {\mathcal H}^+$ and
\[
\text{ran}\bigl(E(\{q^n\}\bigr) = I^\ast {\mathcal H}^+_n = \bigl\{
\text{Per}(f/\sqrt{w})\cdot s_n \,|\, f\in L^2(\Om)\bigr\},
\]
where $s_n$ is the orthonormal Stieltjes--Wigert polynomial of degree
$n$. Note that by the functional equation \eqref{eq:func-eq}, we have
\begin{equation*}
\text{Per}\bigl({f}/{\sqrt{w}}\bigr)(x) =\frac{(Pf)(x)}{\sqrt{w(x)}},
\qquad (Pf)(x) = \sum_{k=-\infty}^\infty
\chi_{(q^{k+1},q^k]}(x)\,  x^{k/2} q^{-k(k+1)/4} f(xq^{-k})
\end{equation*}
and $(Pf)(xq)= \sqrt{x} (Pf)(x)$.  In particular, by taking any
orthonormal basis $\{f_j\}_{j\in\N}$ of $L^2(\Om)$ we obtain from the
orthonormality of $t\mapsto f_j(t)\psi(q^n;t)/N_{q^n}(t)$ in
${\mathcal H}^+$ and the unitarity of $I$ the orthogonality relations
\begin{align}
\label{eq:orthorelforSWpols}
\notag
& \int_0^\infty \text{Per}\bigl({f_i}/{\sqrt{w}}\bigr)(x)\,
\text{Per}\bigl({f_j}/{\sqrt{w}}\bigr)(x)\, s_n(x) s_m(x)\, w(x) dx \\
& \qquad\quad=
\int_{\text{supp}(\mu)} (Pf_i)(x)\,(Pf_j)(x)\, s_n(x) s_m(x) dx =
\de_{n,m}\de_{i,j}.
\end{align}
The special case $i=j$ tells us that the Stieltjes--Wigert polynomials
are orthogonal with respect to any absolutely continuous measure whose
density satisfies the functional equation \eqref{eq:func-eq}. This
result is also obtained in \cite[Prop.~2.1]{ChriJMAA}.

To sum up, we denote by $\text{PPol}\subset L^2(\mu)$ the closure of
the space of functions of the form $\sum f_n p_n\in L^2(\mu)$, with
$f_n$ a $q$-periodic function and $p_n$ a polynomial. It follows that
$\text{PPol}=I^\ast {\mathcal H}^+\subset\text{dom}(L^\ast)$ and
$L^\ast\vert_{\text{PPol}}$ is a bounded linear operator on
$\text{PPol}$ with spectrum $q^{\Z_+}\cup \{0\}$.

We now take a closer look at the spectral decomposition corresponding
to the continuous spectrum of $L^\ast$. For any Borel set $A\subset
(-q^{l-1},-q^l]$ we have $E(A) I^\ast{\mathcal H}^-_r= \{0\}$ unless
$r=l$. Since $E(A)F=I^\ast E(J^\ast| A) IF$ for $F\in L^2(\mu)$ with
compact support, it thus follows that
\[
E(J^\ast|A)(IF)(t)=
\chi_{A \cap \tilde{\Om}_l}(-q^l/t)
\frac{\bigl\langle (IF)(t),\psi(-q^l/t;t)\bigr\rangle_{\ell^2(\Z)}}
{N_{-q^l/t}(t)}\frac{\psi(-q^l/t;t)}{N_{-q^l/t}(t)}.
\]
Calculating $I^\ast$ on ${\mathcal H}^-_l$ gives 
\begin{align*}
\label{eq:IastonHmin}
\notag
I^\ast\biggl(t\mapsto f(t)\frac{\psi(-q^l/t;t)}{N_{-q^l/t}(t)}\biggr)(x) &=
I^\ast\biggl(t\mapsto  \frac{f(t)}{N_{-q^l/t}(t)}
\sum_{k=-\infty}^\infty t^{k/2}
q^{k(k+1)/4} \phi_{-q^l/t}(tq^k) e_k\biggr)(x) \\
&= \sum_{k=-\infty}^\infty \chi_{(q^{k+1},q^k]}(x)
\frac{f(xq^{-k})x^{k/2}}{N_{-q^{l+k}/x}(xq^{-k})}
\frac{ q^{-k(k+1)/4}}{\sqrt{w(x)}} \phi_{-q^{l+k}/x}(x),
\end{align*}
so when $f$ has the form
\begin{equation*}
f(t)=\chi_{A \cap \tilde{\Om}_l}(-q^l/t)
\frac{\bigl\langle (IF)(t), \psi(-q^l/t;t)\bigr\rangle_{\ell^2(\Z)}}
{{N_{-q^l/t}(t)}},
\end{equation*}
we obtain for $G\in L^2(\mu)$ with compact support that
\begin{align*}
\langle E(A)F,G\rangle_{L^2(\mu)} &=
\int_0^\infty \bigl(I^\ast E(J^\ast|A) IF\bigr)(x)\,
\overline{G(x)}\,w(x)dx \\
&=\int_0^\infty \sum_{k=-\infty}^\infty
\chi_{(q^{k+1},q^k]}(x)\frac{\chi_A(-q^{l+k}/{x})}
{N_{-q^{l+k}/x}(xq^{-k})^2}\,
x^{k/2} q^{-k(k+1)/4} \phi_{-q^{l+k}/x}(x) \\
&\qquad\qquad\times
\bigl\langle (IF)(xq^{-k}), \psi(-q^{l+k}/x;xq^{-k})\bigr\rangle_{\ell^2(\Z)}
\overline{G(x)} \sqrt{w(x)} dx.
\end{align*}
Expanding the inner product in the integrand, the integral can be
written as
\begin{align}
\label{eq:spectraldecomp1}
\notag
&\sum_{j,k=-\infty}^\infty \int_{q^{k+1}}^{q^k} \chi_A(-q^{l+k}/{x})
\frac{\phi_{-q^{l+k}/x}(x) \phi_{-q^{l+k}/x}(xq^{j-k})}
{N_{-q^{l+k}/x}(xq^{-k})^2}\, x^j q^{\binom{j+1}{2}-jk}
F(xq^{j-k})\overline{G(x)}w(x) dx \\
\notag
&\quad= \sum_{j,k=-\infty}^\infty (-1)^{j+k}
q^{\binom{j+1}{2}\binom{k+1}{2}+l(j+k)}
\int_A \frac{\phi_\la(-q^{l+j}/\la)\phi_\la(-q^{l+k}/\la)}
{N_\la(-q^l/\la)^2 \, \la^{j+k}} \\
\notag
&\quad\qquad\qquad\qquad\qquad\qquad\qquad\qquad\qquad \times
F(-q^{l+j}/\la) \overline{G(-q^{l+k}/\la)}\,
w(-q^l/\la)\frac{q^{l}}{\la^2}\, d\la \\
\notag
&\quad= \int_A \biggl( \sum_{j=-\infty}^\infty
F(-q^{l+j}/\la)(-q^l/\la)^j q^{\binom{j+1}{2}}
\phi_\la(-q^{l+j}/\la)\biggr)  \\
& \quad\qquad\quad\times
\biggl( \sum_{k=-\infty}^\infty
\overline{G(-q^{l+k}/\la)} (-q^l/\la)^k q^{\binom{k+1}{2}}
\phi_\la(-q^{l+k}/\la)\biggr)
\frac{q^l w(-q^l/\la)}{\la^2 N_\la(-q^l/\la)^2} \, d\la,
\end{align}
using the functional equation \eqref{eq:func-eq}, switched to
$\la=-q^{l+k}/x$. Note that
\[
\sum_{j=-\infty}^\infty F(-q^{l+j}/\la)
(-q^l/\la)^j q^{\binom{j+1}{2}} \phi_\la(-q^{l+j}/\la)=
\frac{(-\la)^l}{q^{\binom{l+1}{2}}}
\sum_{j=-\infty}^\infty F(-q^j/\la)(-\la)^{-j}
q^{\binom{j+1}{2}} \phi_\la(-q^j/\la)
\]
and define
\begin{equation}
\label{eq:defF}
\bigl( {\mathcal F}F\bigr)(\la)=
\sum_{j=-\infty}^\infty F(-q^j/\la)(-\la)^{-j}
q^{\binom{j+1}{2}}\phi_\la(-q^j/\la).
\end{equation}
By means of \eqref{eq:defF} we can write \eqref{eq:spectraldecomp1} as
\begin{align}
\label{eq:spectraldecomp2}
\notag
\langle E(A)F,G\rangle_{L^2(\mu)} &=
\int_A \bigl({\mathcal F}F\bigr)(\la)\overline{\bigl({\mathcal
    F}G\bigr)(\la)}\,
\la^{2l}q^{-l(l+1)}\frac{q^l w(-q^l/\la)}{\la^2 N_\la(-q^l/\la)^2}\,d\la \\
&= \int_A \bigl({\mathcal F}F\bigr)(\la)\overline{\bigl( {\mathcal
      F}G\bigr)(\la)}\,
|\la|^l q^{-l(l+1)/2} \frac{w(-1/\la)}{N_\la(-q^l/\la)^2}\,
\frac{d\la}{\la^2},
\end{align}
using the functional equation \eqref{eq:func-eq} once more. Now define
\begin{equation}
\label{eq:nu}
\nu(\la) = \sum_{l=-\infty}^\infty \chi_{(-q^{l-1}, -q^l]}(\la)
\frac{|\la|^l q^{-l(l+1)/2}}{N_\la(-q^l/\la)^2}
\end{equation}
and use \eqref{eq:spectraldecomp2} to obtain
\begin{equation}\label{eq:spectraldecomp3}
\langle E(A)F,G\rangle_{L^2(\mu)}=
\int_A \bigl( {\mathcal F}F\bigr)(\la)
\overline{\bigl({\mathcal F}G\bigr)(\la)}\,
\nu(\la) w({-1}/{\la}) \, \frac{d\la}{\la^2}
\end{equation}
for an arbitrary Borel set $A\subset (-\infty,0)$. It follows that the
complex measure $\langle E(A)F,G\rangle_{L^2(\mu)}$ is absolutely
continuous with respect to the Lebesgue measure on $(0,\infty)$, and
for any $F, G\in I^\ast {\mathcal H}^-$ we have
\begin{equation}\label{eq:spectraldecomp4}
\langle F,G\rangle_{L^2(\mu)}=\int_{-\infty}^0
\bigl({\mathcal F}F\bigr)(\la)\overline{\bigl({\mathcal F}G\bigr)(\la)}
\,\nu(\la)w({-1}/{\la})\,\frac{d\la}{\la^2}.
\end{equation}

Taking into account the discrete spectrum of $L^\ast$ on the space
$\text{PPol}$  as well, we obtain the following Plancherel type theorem.

\begin{thm}
\label{thm:Plancherel}
Consider an absolutely continuous positive measure $\mu$ on
$(0,\infty)$ with density $w$ satisfying the functional equation
\eqref{eq:func-eq}. Let $\Om=(q,1]\cap\text{\rm supp}(\mu)$ and
suppose that $\{f_i\}_{i=0}^\infty$ is an arbitrary fixed orthonormal
basis of $L^2( \Om)$. For all $F,G\in L^2(\mu)$, we have the
Plancherel equality
\[
\int_0^\infty F(x)\overline{G(x)}w(x)dx =
\sum_{i,n=0}^\infty F_{in} \overline{G_{in}} +
\int_{-\infty}^0 \bigl( {\mathcal F}F\bigr)(\la)
\overline{\bigl( {\mathcal F}G\bigr)(\la)}\,
\nu(\la) w({-1}/{\la}) \, \frac{d\la}{\la^2},
\]
where
\[
F_{in}=\int_0^\infty F(x) \text{\rm Per}
\bigl({f_i}/{\sqrt{w}}\bigr)(x) s_n(x) w(x)\, dx
\]
and ${\mathcal F}$, respectively $\nu$, are defined in \eqref{eq:defF}
and \eqref{eq:nu}.
\end{thm}

We can rewrite the above result in terms of a corresponding transform.
Consider the Hilbert space
\[
{\mathcal K}=\ell^2\bigl(\Z_+\times\Z_+\bigr) \oplus
L^2\Bigl((-\infty,0), \nu(\la)w({-1}/{\la})\frac{d\la}{\la^2}\Bigl)
\]
and define
\begin{equation}
\label{eq:adjointF}
({\mathcal F}^\ast g)(x) = \sum_{i,n=0}^\infty g_{in}\,
\text{\rm Per}\bigl({f_i}/{\sqrt{w}}\bigr)(x) s_n(x)
+\sum_{j=-\infty}^\infty g(-{q^j}/{x})\,
\phi_{-q^j/x}(x) \nu(-{q^j}/{x}), \quad x>0
\end{equation}
for compactly supported functions $g\in {\mathcal K}$.
If we consider ${\mathcal F}$ as defined in \eqref{eq:defF} as
${\mathcal F} \colon I^\ast {\mathcal H}^- \to L^2\Bigl((-\infty,0),
\nu(\la)w({-1}/{\la})\frac{d\la}{\la^2}\Bigl)$ and extend it to an
operator ${\mathcal F} \colon L^2(\mu)\to {\mathcal K}$ by defining
${\mathcal F}\colon I^\ast{\mathcal H}^+\to
\ell^2\bigl(\Z_+\times\Z_+\bigr)$ by ${\mathcal F}F= \{
F_{in}\}_{i,n\in\Z_+}$ with $F_{in}$ as in Theorem
\ref{thm:Plancherel}, then we have the following result.

\begin{cor}
\label{cor:thm:Plancherel}
${\mathcal F}
\colon L^2(\mu) \to {\mathcal K}$ is a unitary isomorphism with
adjoint given by \eqref{eq:adjointF}.
\end{cor}

\def\cprime{$'$}
\providecommand{\bysame}{\leavevmode\hbox to3em{\hrulefill}\thinspace}
\providecommand{\MR}{\relax\ifhmode\unskip\space\fi MR }
\providecommand{\MRhref}[2]{%
  \href{http://www.ams.org/mathscinet-getitem?mr=#1}{#2}
}
\providecommand{\href}[2]{#2}


\begin{thebibliography}{10}

\bibitem{Akhi}
Naum~Ilyich Akhiezer, \emph{The classical moment problem and some related
  questions in analysis}, Translated by N. Kemmer, Hafner Publishing Co., New
  York, 1965.

\bibitem{Nodarse+Medem}
Renato {\'A}lvarez-Nodarse and Juan~Carlos Medem, \emph{{$q$}-classical
  polynomials and the {$q$}-{A}skey and {N}ikiforov-{U}varov tableaus}, J.
  Comput. Appl. Math. \textbf{135} (2001), no.~2, 197--223.

\bibitem{Andrews}
George~E. Andrews, \emph{Ramanujan's ``lost'' notebook. {VIII}: {T}he entire
  {R}ogers-{R}amanujan function}, Adv. Math. \textbf{191} (2005), no.~2,
  393--407.

\bibitem{Bere}
Yurij~M. Berezans{\cprime}ki{\u\i}, \emph{Expansions in eigenfunctions of
  selfadjoint operators}, Translated from the Russian by R. Bolstein, J. M.
  Danskin, J. Rovnyak and L. Shulman. Translations of Mathematical Monographs,
  Vol. 17, American Mathematical Society, Providence, R.I., 1968.

\bibitem{Berg}
Christian Berg, \emph{From discrete to absolutely continuous solutions of
  indeterminate moment problems}, Arab J. Math. Sci. \textbf{4} (1998), no.~2,
  1--18.

\bibitem{Berg-geo}
\bysame, \emph{On some indeterminate moment problems for measures on a
  geometric progression}, J. Comput. Appl. Math. \textbf{99} (1998), no.~1,
  67--75.

\bibitem{Carlitz}
Leonard Carlitz, \emph{Fibonacci notes. {IV}. {$q$}-{F}ibonacci polynomials},
  Fibonacci Quart. \textbf{13} (1975), 97--102.

\bibitem{Chihara}
Theodore~Seio Chihara, \emph{A characterization and a class of distribution
  functions for the {S}tieltjes-{W}igert polynomials}, Canad. Math. Bull.
  \textbf{13} (1970), 529--532.

\bibitem{ChriCA}
Jacob~Stordal Christiansen, \emph{The moment problem associated with the
  {$q$}-{L}aguerre polynomials}, Constr. Approx. \textbf{19} (2003), no.~1,
  1--22.

\bibitem{ChriJMAA}
\bysame, \emph{The moment problem associated with the {S}tieltjes-{W}igert
  polynomials}, J. Math. Anal. Appl. \textbf{277} (2003), no.~1, 218--245.

\bibitem{ChriPhD}
\bysame, \emph{{I}ndeterminate moment problems within the {A}skey-scheme},
  Ph.D. thesis, University of Copenhagen (2004),
  http://www.math.ku.dk$\sim$stordal/thesis.pdf.

\bibitem{CiccKK}
Nicola Ciccoli, Erik Koelink, and Tom~H. Koornwinder, \emph{{$q$}-{L}aguerre
  polynomials and big {$q$}-{B}essel functions and their orthogonality
  relations}, Methods Appl. Anal. \textbf{6} (1999), no.~1, 109--127.

\bibitem{Dixm}
Jacques Dixmier, \emph{von {N}eumann algebras}, North-Holland Mathematical
  Library, vol.~27, North-Holland Publishing Co., Amsterdam, 1981, Translated
  from the second French edition by F. Jellett.

\bibitem{DunfS}
Nelson Dunford and Jacob~T. Schwartz, \emph{Linear operators. {P}art {II}},
  Wiley Classics Library, John Wiley \& Sons Inc., New York, 1988, Spectral
  theory. Selfadjoint operators in Hilbert space, Reprint of the 1963 original,
  A Wiley-Interscience Publication.

\bibitem{GaspR}
George Gasper and Mizan Rahman, \emph{Basic hypergeometric series, 2nd ed.},
  Encyclopedia of Mathematics and its Applications, vol.~96, Cambridge
  University Press, Cambridge, 2004.

\bibitem{Hayman}
Walter~K. Hayman, \emph{On the zeros of a $q$-{B}essel function}, to appear in
  Contemp. Math.

\bibitem{Mourad}
Mourad E.~H. Ismail, \emph{On {J}ackson's third $q$-{B}essel function and
  $q$-exponentials}, Preprint (2000).

\bibitem{K&S}
Roelof Koekoek and Ren\'e~F. Swarttouw, \emph{The {A}skey-scheme of
  hypergeometric orthogonal polynomials and its $q$-analogue}, Tech. Report no.
  98-17, TU--Delft (1998).

\bibitem{KoelLaredo}
Erik Koelink, \emph{Spectral theory and special functions}, Laredo Lectures on
  Orthogonal Polynomials and Special Functions, Adv. Theory Spec. Funct.
  Orthogonal Polynomials, Nova Sci. Publ., Hauppauge, NY, 2004, pp.~45--84.

\bibitem{KoelJMAA}
H.~T. Koelink, \emph{Hansen-{L}ommel orthogonality relations for {J}ackson's
  {$q$}-{B}essel functions}, J. Math. Anal. Appl. \textbf{175} (1993), no.~2,
  425--437.

\bibitem{Marcellan+Medem}
Francisco Marcell{\'a}n and Juan~Carlos Medem, \emph{{$q$}-classical orthogonal
  polynomials: a very classical approach}, Electron. Trans. Numer. Anal.
  \textbf{9} (1999), 112--127 (electronic), Orthogonal polynomials: numerical
  and symbolic algorithms (Legan\'es, 1998).

\bibitem{MassR}
David~R. Masson and Joe Repka, \emph{Spectral theory of {J}acobi matrices in
  {$l\sp 2({\bf Z})$} and the {${\rm su}(1,1)$} {L}ie algebra}, SIAM J. Math.
  Anal. \textbf{22} (1991), no.~4, 1131--1146.

\bibitem{ReedS4}
Michael Reed and Barry Simon, \emph{Methods of modern mathematical physics.
  {IV}. {A}nalysis of operators}, Academic Press [Harcourt Brace Jovanovich
  Publishers], New York, 1978.

\end{thebibliography}
\end{document}